\documentclass[12pt, reqno]{amsart}

\usepackage[margin=.9in, a4paper, foot=.3in]{geometry}

\usepackage{graphicx, color}

\definecolor{darkblue}{rgb}{0,0,.5}
\definecolor{darkred}{rgb}{.5,0,0}
\definecolor{darkgreen}{rgb}{0,0.5,0}

\usepackage{hyperref}
\hypersetup {
    pdfstartview=FitH,
    colorlinks,
    citecolor=darkgreen,
    linkcolor=darkred,
    urlcolor=darkblue
}

\usepackage[noBBpl,slantedGreek]{mathpazo}

\usepackage[mathcal]{euscript}

\usepackage{amssymb}

\usepackage{bm}

\newcommand {\bbB}{\mathbb B}
\newcommand {\bbC}{\mathbb C}
\newcommand {\bbD}{\mathbb D}
\newcommand {\bbE}{\mathbb E}
\newcommand {\bbK}{\mathbb K}
\newcommand {\bbM}{\mathbb M}
\newcommand {\bbN}{\mathbb N}
\newcommand {\bbR}{\mathbb R}

\newcommand {\bbU}{\mathbb U}
\newcommand {\bbV}{\mathbb V}
\newcommand {\bbW}{\mathbb W}
\newcommand {\bbZ}{\mathbb Z}

\newcommand {\mbar}[3]{\hskip #2 \overline{\hskip -#2 #1 \hskip -#3} \hskip #3}

\newcommand {\oD}{\mbar{D}{.15em}{.1em}}
\newcommand {\oM}{\mbar{M}{.2em}{.1em}}
\newcommand {\oU}{\mbar{U}{.05em}{.1em}}
\newcommand {\oV}{\mbar{V}{.05em}{.1em}}
\newcommand {\ovarphi}{\mbar{\varphi}{.1em}{.05em}}
\newcommand {\oW}{\mbar{W}{.05em}{.05em}}

\newcommand {\obbD}{\mbar{\mathbb D}{0em}{.1em}}
\newcommand {\obbM}{\mbar{\mathbb M}{.05em}{.1em}}
\newcommand {\obbU}{\mbar{\mathbb U}{0em}{.06em}}
\newcommand {\obbV}{\mbar{\mathbb V}{0em}{.05em}}
\newcommand {\obbW}{\overline {\mathbb W}}

\newcommand {\calB}{\mathcal B}

\newcommand {\calK}{\mathcal K}
\newcommand {\calL}{\mathcal L}
\newcommand {\calR}{\mathcal R}

\newcommand {\gothg}{\mathfrak g}
\newcommand {\gothgl}{\mathfrak{gl}}
\newcommand {\gothh}{\mathfrak h}

\newcommand {\gothsl}{\mathfrak{sl}}

\newcommand {\rme}{\mathrm e}

\newcommand {\tlsliii}{\widetilde{\mathcal L}(\mathfrak{sl}_3)}
\newcommand {\uqbm}{\mathrm U_q(\mathfrak b_-)}
\newcommand {\uqbp}{\mathrm U_q(\mathfrak b_+)}
\newcommand {\uqnm}{\mathrm U_q(\mathfrak n_-)}
\newcommand {\uqnp}{\mathrm U_q(\mathfrak n_+)}
\newcommand {\uqgliii}{\mathrm U_q(\mathfrak{gl}_3)}
\newcommand {\uqsliii}{\mathrm U_q(\mathfrak{sl}_3)}
\newcommand {\uqlsliii}{\mathrm U_q(\mathcal L(\mathfrak{sl}_3))}
\newcommand {\uqtlsliii}{\mathrm U_q(\widetilde{\mathcal L}(\mathfrak{sl}_3))}

\makeatletter
\def\section{\@startsection{section}{1}%
  \z@{-.7\linespacing\@plus -\linespacing}{.5\linespacing}%
  {\normalfont\scshape\centering}}
\def\subsection{\@startsection{subsection}{2}%
  \z@{-.5\linespacing\@plus -.7\linespacing}{.5em}%
  {\normalfont\bfseries\mathversion{bold}}}
\makeatother

\allowdisplaybreaks

\numberwithin{equation}{section}

\DeclareMathOperator {\End}{End}
\DeclareMathOperator {\id}{id}
\DeclareMathOperator {\Mat}{Mat}

\title{Monodromy operators for higher rank}

\author[A. V. Razumov]{A. V. Razumov}

\address{Institute for High Energy Physics, 142281 Protvino, Moscow region, Russia}

\email{Alexander.Razumov@ihep.ru}

\begin{document}

\addtolength {\jot}{3pt}

\begin{abstract}
We find the explicit form of the basic monodromy operators for the case of the quantum group $\mathrm U_q(\mathcal L(\mathfrak{sl}_3))$. Expressions for the quantum Casimir elements of the quantum groups $\mathrm U_q(\mathfrak{gl}_3)$ and $\mathrm U_q(\mathfrak{sl}_3)$ are obtained as a by-product.
\end{abstract}

\maketitle

\hbox to \hsize{\hfil \em To the memory of Yuri~Stroganov.}

\section{Introduction}

The most productive, although not comprehensive, approach to the theory of quantum integrable systems is based on the concept of a quantum group introduced by Drinfeld and Jimbo \cite{Dri87, Jim85}. In this approach, all the objects describing the model and related to its integrability originate from the universal $R$-matrix. For the first time, this was consistently realised by Bazhanov, Lukyanov and Zamolodchikov \cite{BazLukZam96, BazLukZam97, BazLukZam99}, see also \cite{AntFei97, BooGoeKluNirRaz12a, BooGoeKluNirRaz13}. The method was used to obtain an explicit form of $R$-matrices for certain representations of the quantum groups $\mathrm U_q(\calL(\gothsl_2))$ \cite{KhoTol92, LevSoiStu93, ZhaGou94, BraGouZhaDel94, BraGouZha95, BooGoeKluNirRaz10}, $\mathrm U_q(\calL(\gothsl_3))$ \cite{ZhaGou94, BraGouZhaDel94, BraGouZha95, BooGoeKluNirRaz10} and $\mathrm U_q(\calL(\gothsl_3, \mu))$ \cite{KhoTol92, BooGoeKluNirRaz11}, where $\mu$ is the standard diagram automorphism of the Lie algebra $\gothsl_3$ of
order 2. An example of a quantum supergroup was considered in the paper \cite{BazTsu08}. It appears that one can also find the form of monodromy operators, transfer matrices, $L$-operators, and $Q$-operators \cite{BazTsu08, BooGoeKluNirRaz10, BooGoeKluNirRaz12a}.

The universal $R$-matrix is an element of the tensor product of two copies of the quantum group under consideration. A monodromy operator for a discrete quantum integrable system is constructed by a choice of two finite-dimensional representations of the quantum group.\footnote{In fact, the case when one of the representations is infinite-dimensional is also interesting.} For the case of the quantum group $\mathrm U_q(\calL(\gothsl_n))$ the usual way to construct finite-dimensional representations is to use the homomorphism from $\mathrm U_q(\calL(\gothsl_n))$ to $\mathrm U_q(\gothgl_n)$ proposed by Jimbo \cite{Jim86a}, see, for example \cite{ChaPre91, ChaPre94}. It is convenient for applications to fix the representation only for the second factor of the tensor product, and to use for the first factor the Jimbo's homomorphism. Here the monodromy operator can be considered as a matrix with entries in the quantum group $\mathrm U_q(\gothgl_n)$. For the case of the quantum group $\mathrm U_q(\calL(\gothsl_2))$ 
such monodromy operator was obtained in the paper \cite{BooGoeKluNirRaz12a}. In the present paper we consider the case of $\mathrm U_q(\calL(\gothsl_3))$. It is worth to note that the general form of the monodromy operator in question up to a factor belonging to the centre of~$\mathrm U_q(\gothgl_n)$ was found by Jimbo \cite{Jim86a}. Our main goal is to find this factor. As a by-product we obtain expressions for the quantum Casimir elements of the quantum groups $\mathrm U_q(\gothgl_3)$ and $\mathrm U_q(\gothsl_3)$.

Below $\bbZ$ denotes the set of integers, $\bbZ_+$ the set of non-negative integers, and $\bbN$ the set of positive integers. Depending on the context, the symbol `$1$' means the number one, the unit of an algebra or the unit matrix. We use the notation
\begin{equation*}
\kappa_q = q - q^{-1},
\end{equation*}
so that the $q$-deformation of a number $\nu \in \bbC$ is
\begin{equation*}
[\nu]_q = \frac{q^\nu - q^{-\nu}}{q - q^{-1}} = \kappa_q^{-1} (q^\nu - q^{-\nu}).
\end{equation*}
When necessary we identify a linear operator with its matrix with respect to a basis. At last, following a physical tradition, we sometimes call elements of an algebra operators.

\section{\texorpdfstring{Quantum groups $\uqgliii$ and $\uqsliii$}{Quantum groups Uq(gl3) and Uq(sl3)}}

Depending on the sense of the deformation parameter $q$, there are at least three definitions of a quantum group. According to the first definition, $q = \exp \hbar$, where $\hbar$ is an indeterminate, according to the second one, $q$ is indeterminate, and according to the third one, $q = \exp \hbar$, where $\hbar$ is a complex number. In the first case a quantum group is a $\bbC[[\hbar]]$-algebra, in the second case a $\bbC(q)$-algebra, and in the third case it is just a complex algebra. Usually one uses monodromy operators to construct transfer operators with the help of some trace operation. It seems that to this end it is convenient to use the third definition of a quantum group. Therefore, we define the quantum group as a $\bbC$-algebra, see, for example, the books \cite{JimMiw85, EtiFreKir98}.

Denote by $\gothg$ the standard Cartan subalgebra of the Lie algebra $\gothgl_3$ and by $G_i$, $i = 1, 2, 3$,  the standard Cartan generators. The root system of $\gothgl_3$ relative to $\gothg$ is generated by the simple roots $\alpha_i \in \gothg^*$, $i = 1, 2$, given by the relations
\begin{equation}
\alpha_j(G_i) = c_{ij}, \label{alphah}
\end{equation}
where
\begin{equation*}
(c_{ij}) = \left( \begin{array}{rr}
1 & 0 \\
-1 & 1 \\
0 & -1
\end{array} \right).
\end{equation*}
The Lie algebra $\gothsl_3$ is a subalgebra of $\gothgl_3$, and the standard Cartan subalgebra $\gothh$ of $\gothsl_3$ is a subalgebra of $\gothg$. Here the standard Cartan generators $H_i$, $i = 1, 2$, of $\gothsl_3$ are
\begin{equation}
H_1 = G_1 - G_2, \qquad H_2 = G_2 - G_3, \label{hk}
\end{equation}
and we have
\begin{equation*}
\alpha_j(H_i) = a_{ij},
\end{equation*}
where
\begin{equation}
(a_{ij}) = \left( \begin{array}{rr}
2 & -1 \\
-1 & 2
\end{array} \right) \label{cm}
\end{equation}
is the Cartan matrix of $\gothsl_3$.

Let $\hbar$ be a complex number and $q = \exp \hbar$. The quantum group $\uqgliii$ is a unital associative $\bbC$-algebra generated by the elements $E_i$, $F_i$, $i = 1, 2$, and $q^X$, $X \in \gothg$, with the relations\footnote{Here and below we assume that $q^2 \ne 1$.}
\begin{gather}
q^0 = 1, \qquad q^{X_1} q^{X_2} = q^{X_1 + X_2}, \label{xx} \\
q^X E_i q^{-X} = q^{\alpha_i(X)} E_i, \qquad q^X F_i q^{-X} = q^{- \alpha_i(X)} F_i, \label{xexf} \\
[E_i, F_j] = \delta_{ij} \frac{q^{H_i} - q^{-H_i}}{q - q^{-1}} \label{ef}
\end{gather}
satisfied for any $i$ and $j$, and the Serre relations
\begin{equation}
E_i^2 E_j^{} - [2]_q E_i^{} E_j^{} E_i^{} + E_j^{} E_i^2 = 0, \qquad
F_i^2 F_j^{} - [2]_q F_i^{} F_j^{} F_i^{} + F_j^{} F_i^2 = 0 \label{sr}
\end{equation}
satisfied for distinct $i$ and $j$. Note that $q^X$ is just a convenient notation. There are no elements of $\uqgliii$ corresponding to the elements of $\gothg$. Below for any $X \in \gothg$ and $\nu \in \bbC$ we use the notation
\begin{equation*}
q^{X + \nu} = q^\nu q^X.
\end{equation*}

Looking at (\ref{xexf}) one can say that the generators $e_i$ and $f_i$ are related to the roots $\alpha_i$ and $-\alpha_i$ respectively. We define the elements related to the roots $\alpha_1 + \alpha_2$ and $-(\alpha_1 + \alpha_2)$ as
\begin{equation}
E_3 = E_1 E_2 - q^{-1} E_2 E_1, \qquad F_3 = F_2 F_1 - q F_1 F_2. \label{e3f3}
\end{equation}
With respect to the properly defined coproduct, counit and antipode the quantum group $\uqgliii$ is a Hopf algebra.

There is a useful set of automorphisms of $\uqsliii$ defined as
\begin{gather}
E_i \to \nu_i^{\mathstrut} E^{\mathstrut}_i \, q^{\sum_{j = 1}^3 \nu_{ij} G_j}, \qquad F_i \to \nu_i^{-1} q^{- \sum_{j = 1}^3 \nu_{ij} G_j} F^{\mathstrut}_i, \qquad q^X \to q^X, \label{aut}
\end{gather}
where $\nu_i$ are arbitrary nonzero complex numbers, and the complex numbers $\nu_{ij}$ satisfy the relations
\begin{equation*}
\sum_{k = 1}^3 \nu_{ik} c_{kj} = \sum_{k = 1}^3 \nu_{jk} c_{ki}.
\end{equation*}
It is easy to see that all these relations are equivalent to only one equality
\begin{equation*}
\nu_{12} - \nu_{13} = \nu_{21} - \nu_{22}.
\end{equation*}

The quantum group $\uqsliii$ can be defined as the subalgebra of $\uqgliii$ generated by $E_i$, $F_i$, $i = 1, 2$, and $q^X$, $X \in \gothh$. Actually, to construct monodromy operators one can use the corresponding homomorphism from $\uqlsliii$ to $\uqsliii$, which is a modification of the Jimbo's homomorphism. However, the expressions obtained in this way are more complicated, and we will work with the homomorphism to the quantum group $\uqgliii$.In fact, the expressions for the case of $\uqsliii$ can be easily deduced from the expressions for the case of $\uqgliii$, see appendix.

\section{\texorpdfstring{Quantum group $\uqlsliii$}{Quantum group Uq(L(sl3))}}

\subsection{Definition}

We start our consideration of the quantum group $\uqlsliii$ with the quantum group $\uqtlsliii$. Remind that the Cartan subalgebra of $\tlsliii$ is
\begin{equation*}
\widetilde \gothh = \gothh \oplus \bbC c,
\end{equation*}
where $\gothh = \bbC H_1 \oplus \bbC H_2$ is the standard Cartan subalgebra of the Lie algebra $\gothsl_3$ and $c$ is the central element \cite{Kac90}. Define the Cartan elements
\begin{equation*}
h_0 = c - H_1 - H_2, \qquad h_1 = H_1, \qquad h_2 = H_2,
\end{equation*}
so that one has
\begin{equation*}
\widetilde \gothh = \bbC h_0 \oplus \bbC h_1 \oplus \bbC h_2.
\end{equation*}
The simple roots $\alpha_i \in \widetilde{\gothh}^*$, $i = 0, 1, 2$, are given by the equality
\begin{equation*}
\alpha_j(h_i) = \tilde a_{ij},
\end{equation*}
where
\begin{equation*}
(\tilde a_{ij}) = \left(\begin{array}{rrr}
2 & -1 & -1 \\
-1 & 2 & -1 \\
-1 & -1 & 2
\end{array} \right)
\end{equation*}
is the Cartan matrix of the Lie algebra $\tlsliii$.

As before, let $\hbar$ be a complex number and $q = \exp \hbar$. The quantum group $\uqtlsliii$  is a unital associative $\bbC$-algebra generated by the elements $e_i$, $f_i$, $i = 0, 1, 2$, and $q^x$, $x
\in \widetilde \gothh$, with the relations
\begin{gather}
q^0 = 1, \qquad q^{x_1} q^{x_2} = q^{x_1 + x_2}, \label{lxx} \\
q^x e_i q^{-x} = q^{\alpha_i(x)} e_i, \qquad q^x f_i q^{-x} = q^{-\alpha_i(x)} f_i, \label{lxexf} \\
[e_i, f_j] = \delta_{ij} \, \frac{q^{h_i} - q^{-h_i}}{q - q^{-1}} \label{lef}
\end{gather}
satisfied for all $i$ and $j$, and the Serre relations
\begin{equation}
e_i^2 e_j^{\mathstrut} - [2]_q  e_i^{\mathstrut} e_j^{\mathstrut} e_i^{\mathstrut}
+ e_j^{\mathstrut} e_i^2 = 0, \qquad f_i^2 f_j^{\mathstrut} - [2]_q  f_i^{\mathstrut} f_j^{\mathstrut} f_i^{\mathstrut} + f_j^{\mathstrut} f_i^2 = 0 \label{lsr}
\end{equation}
satisfied for distinct $i$ and $j$.

The quantum group $\uqtlsliii$ is a Hopf algebra with the comultiplication $\Delta$ defined by the
relations
\begin{gather}
\Delta(q^x) = q^x \otimes q^x, \label{dqx} \\
\Delta(e_i) = e_i \otimes 1 + q^{-h_i} \otimes e_i, \qquad \Delta(f_i) = f_i \otimes q^{h_i} + 1
\otimes f_i, \label{defi}
\end{gather}
and with the correspondingly defined counit and antipode.

The quantum group $\uqlsliii$ can be defined as the quotient algebra of $\uqtlsliii$ by the two-sided
ideal generated by the elements of the form $q^{\nu c} - 1$, $\nu \in \bbC^\times$. In terms of
generators and relations the quantum group $\uqlsliii$ is a $\bbC$-algebra generated by the elements
$e_i$, $f_i$, $i = 0, 1, 2$, and $q^x$, $x \in \widetilde{\gothh}$, with relations
(\ref{lxx})--(\ref{lsr}) and additional relations
\begin{equation}
q^{\nu(h_0 + h_1 + h_2)} = q^{\nu c} = 1, \qquad \nu \in \bbC^\times.\label{qh0h1}
\end{equation}
It is a Hopf algebra with the comultiplication defined by (\ref{dqx}),
(\ref{defi}) and with the correspondingly defined counit and antipode. One of the reasons to use
the quantum group $\uqlsliii$ instead of $\uqtlsliii$ is that in the case of $\uqtlsliii$ we have no expression for the universal $R$-matrix.

\subsection{\texorpdfstring{Universal $R$-matrix}{Universal R-matrix}} \label{s:urm}

As any Hopf algebra the quantum group $\uqlsliii$ has another comultiplication called the opposite comultiplication. It is given by the equalities
\begin{gather}
\Delta^{\mathrm{op}}(q^x) = q^x \otimes q^x, \label{doqx} \\
\Delta^{\mathrm{op}}(e_i) = e_i \otimes q^{-h_i} + 1 \otimes e_i, \qquad \Delta^{\mathrm{op}}(f_i) = f_i \otimes 1 + q^{h_i} \otimes f_i. \label{doefi}
\end{gather}
When the quantum group $\uqlsliii$ is defined as a $\bbC[[\hbar]]$-algebra it is quasitriangular. It means that there exists an invertible element $\calR \in \uqlsliii \otimes \uqlsliii$ such that
\begin{equation}
\Delta^{\mathrm{op}}(a) = \calR \, \Delta(a) \, \calR^{-1} \label{urm}
\end{equation}
for all $a \in \uqlsliii$, and\footnote{For the explanation of the notation and further information see, for example, the book \cite{ChaPre94} or the papers \cite{BooGoeKluNirRaz10, BooGoeKluNirRaz12a}.}
\begin{equation}
(\Delta \otimes \id) (\calR) = \calR^{13} \calR^{23}, \qquad (\id \otimes \Delta) (\calR) = \calR^{13} \calR^{12}. \label{didr}
\end{equation}
The most important property of the universal $R$-matrix is the equality
\begin{equation}
\calR^{12} \calR^{13} \calR^{23} = \calR^{23} \calR^{13} \calR^{12} \label{rrr}
\end{equation}
called the Yang--Baxter equation for the universal $R$-matrix.

The expression for the universal $R$-matrix of $\uqlsliii$ considered as a $\bbC[[\hbar]]$-al\-geb\-ra can be constructed using the procedure proposed by Khoroshkin and Tolstoy~\cite{TolKho92}. Note that here the universal $R$-matrix is an element of $\uqbp \otimes \uqbm$, where $\uqbp$ is the Borel subalgebra of $\uqlsliii$ generated by $e_i$, $i = 0, 1, 2$, and $q^x$, $x \in \widetilde{\gothh}$, and $\uqbm$ is the Borel subalgebra of $\uqlsliii$ generated by $f_i$, $i = 0, 1, 2$, and $q^x$, $x \in \widetilde{\gothh}$.

In fact, one can use the expression for the universal $R$-matrix from the paper \cite{TolKho92} also for the case of a quantum group defined as a $\bbC$-algebras having in mind that in this case the quantum group is quasitriangular only in some restricted sense. Namely, all the relations involving the universal $R$-matrix should be considered as valid only for the weight representations, see in this respect the paper \cite{Tan92}, the book \cite[p. 327]{ChaPre94}, and the discussion below.

A representation $\rho$ of $\uqlsliii$ on the vector space $V$ is a weight representation if
\begin{equation*}
V = \bigoplus_{\lambda \in \widetilde \gothh^*} V_\lambda,
\end{equation*}
where
\begin{equation*}
V_\lambda = \{v \in V \mid q^x v = q^{\lambda(x)} v \mbox{ for any } x \in \widetilde \gothh \}.
\end{equation*}
Note that the element $\lambda$ in the definition of the weight subspace $V_\lambda$ is defined uniquely. Therefore, for a given $x \in \widetilde \gothh$ one can define the operator acting on $v \in V_\lambda$ as the multiplication by $\lambda(x)$. It is natural to denote this operator by $\rho(x)$.

Let $\rho_1$ and $\rho_2$ be two weight representations of $\uqlsliii$ on the vector spaces $V_1$ and $V_2$ with the weight decompositions
\begin{equation*}
V_1 = \bigoplus_{\lambda \in \widetilde \gothh^*} (V_1)_\lambda, \qquad V_2 = \bigoplus_{\lambda \in \widetilde \gothh^*} (V_2)_\lambda.
\end{equation*}
In the tensor product $V_1 \otimes V_2$ the role of the universal $R$-matrix is played by the operator
\begin{equation}
\calR_{\rho_1, \, \rho_2} = (\rho_1 \otimes \rho_2)(\calB) \, \calK_{\rho_1, \, \rho_2}. \label{rpipi}
\end{equation}
Here $\calB$ is an element of $\uqnp \otimes \uqnm$, where $\uqnp$ and $\uqnm$ are the subalgebras of $\uqlsliii$ generated by $e_i$, $i = 0, 1, 2$, and  $f_i$, $i = 0, 1, 2$, respectively. The operator $\calK_{\rho_1, \, \rho_2}$ acts on a vector $v \in (V_1)_{\lambda_1} \otimes (V_2)_{\lambda_2}$ in accordance with the equality
\begin{equation}
\calK_{\rho_1, \, \rho_2} \, v = q^{\sum_{i, j = 1}^2 b_{i j} \, \lambda_1(h_i) \lambda_2(h_j)} \, v, \label{kpipi}
\end{equation}
where
\begin{equation*}
(b_{ij}) = \frac{1}{3} \left( \begin{array}{cc}
2 & 1 \\
1 & 2
\end{array} \right)
\end{equation*}
is the inverse matrix of the Cartan matrix (\ref{cm}) of the Lie algebra $\gothsl_3$. It is clear that
\begin{equation*}
\calK_{\rho_1, \, \rho_2} = q^{\sum_{i, j = 1}^2 b_{ij} \, \rho_1(h_i) \otimes \rho_2(h_j)},
\end{equation*}
and, slightly abusing notation, we write
\begin{equation*}
\calK_{\rho_1, \, \rho_2} = (\rho_1 \otimes \rho_2) \Bigl( q^{\sum_{i, j = 1}^2 b_{ij} \, h_i \otimes h_j} \Bigr) = (\rho_1 \otimes \rho_2)(\calK).
\end{equation*}
Now one can demonstrate that
\begin{equation*}
(\rho_1 \otimes \rho_2)(\Delta^{\mathrm{op}}(a)) = \calR^{\mathstrut}_{\rho_1, \, \rho_2} (\rho_1 \otimes \rho_2)(\Delta(a)) \calR^{-1}_{\rho_1, \, \rho_2}
\end{equation*}
for all $a \in \uqlsliii$. There are similar substitutes for relations (\ref{didr}) and for the Yang--Baxter equation (\ref{rrr}).

In fact, fixing representations $\rho_1$ and $\rho_2$ we obtain objects describing integrable models and helping to analyse them, see, for example, \cite{BooGoeKluNirRaz12a, BooGoeKluNirRaz13}. In particular, if the representations $\rho_1$ and $\rho_2$ coincide with some finite-dimensional representation we have an $R$-operators.  If these representations are finite-dimensional and different we deal with a monodromy operator. For the case of the quantum group $\uqlsliii$ weight representations are usually generated from representations of the quantum group $\uqsliii$ by the Jimbo's homomorphism $\varphi$, see section \ref{s:rm}. Here it is convenient to fix only one representation, say $\rho_2$, and use instead the representation $\rho_1$ the homomorphism $\varphi$. Denote $\rho_2$ by $\rho$ and assume that it is a representation on a vector space $V$. Let $\{e_a\}$ be a basis of $V$ consisting of weight vectors, $\lambda_a$ an element of $\widetilde \gothh^*$ corresponding to the vector $e_a$, and $p_a$
the projection on $e_a$. Now, the role of the universal $R$-matrix is played by the element
\begin{equation}
\calR_{\varphi, \, \rho} = (\varphi \otimes \rho)(\calB) \, \calK_{\varphi, \rho} \label{rphirho}
\end{equation}
of $\uqlsliii \otimes \End(V)$, where
\begin{equation}
\calK_{\varphi, \, \rho} = \sum_a \varphi \Bigl( q^{\sum_{i, j = 1}^2 h_i \, b_{ij} \,\lambda_a(h_j)} \Bigr) \otimes p_a.  \label{kphipi}
\end{equation}
Again, slightly abusing notation, we write
\begin{equation*}
\calK_{\varphi, \, \rho} = (\varphi \otimes \rho) \Bigl( q^{\sum_{i, j = 1}^2 b_{ij} \, h_i \otimes h_j} \Bigr) = (\varphi \otimes \rho)(\calK).
\end{equation*}
It is clear that in the case where $\varphi$ is a weight representation of $\uqlsliii$ this definition is consistent with the definition (\ref{kpipi}).

To describe the structure of the element $\calB$, entering (\ref{rpipi}) and (\ref{rphirho}), we have to define root vectors corresponding to the roots of $\tlsliii$. We say that $a \in \uqtlsliii$ is a root vector corresponding to a root $\gamma$ of $\tlsliii$ if
\begin{equation*}
q^x a \, q^{-x} = q^{\gamma(x)} a
\end{equation*}
for all $x \in \widetilde \gothh$. It is customary to denote
\begin{equation*}
\delta = \alpha_0 + \alpha_1 + \alpha_2, \qquad \alpha = \alpha_1, \qquad \beta = \alpha_2.
\end{equation*}
Let $\triangle_+ = \{\alpha, \beta, \alpha + \beta\}$ denote the system of positive roots of $\gothsl_3$. Then the system of positive roots of $\tlsliii$ is \cite{Kac90}
\begin{multline*}
\widetilde \triangle_+ = \{ \gamma + k \delta \mid \gamma \in \triangle_+, \, \, k \in \bbZ_+ \} \\* \cup \{ k \delta \mid k \in \bbN \} \cup \{ (\delta - \gamma) + k \delta \mid \gamma \in \triangle_+, \, \, k \in \bbZ_+ \}.
\end{multline*}
The full system of roots $\widetilde \triangle$ is the union of the systems of positive and negative roots, $\widetilde \triangle = \widetilde \triangle_+ \cup (- \widetilde \triangle_+)$.

We denote the root vector corresponding to a positive root $\gamma$ by $e_\gamma$, and the root vector corresponding to a negative root $\gamma$ by $f_{-\gamma}$. The evident choice for the root vectors corresponding to the simple roots is
\begin{equation*}
e_{\delta - \alpha - \beta} = e_0, \qquad e_\alpha = e_1, \qquad e_\beta = e_2,
\end{equation*}
and for the simple negative roots is
\begin{equation*}
f_{\delta - \alpha - \beta} = f_0, \qquad f_\alpha = f_1, \qquad f_\beta = f_2.
\end{equation*}
We define the root vectors corresponding to the roots $\pm(\alpha + \beta)$ as
\begin{equation*}
e_{\alpha + \beta} = e_\alpha \, e_\beta - q^{-1} e_\beta \, e_\alpha, \qquad f_{\alpha + \beta} = f_\beta \, f_\alpha - q \, f_\alpha \, f_\beta.
\end{equation*}
and the root vectors corresponding to the roots $\pm (\delta - \gamma)$, $\gamma \in \triangle_+$, as
\begin{align*}
& e_{\delta - \alpha} = e_\beta \, e_{\delta - \alpha - \beta} - q^{-1} e_{\delta - \alpha - \beta} \, e_\beta, && e_{\delta - \beta} = e_\alpha \, e_{\delta - \alpha - \beta} - q^{-1} e_{\delta - \alpha - \beta} \, e_\alpha, \\
& f_{\delta - \alpha} = f_{\delta - \alpha - \beta} \, f_\beta - q \, f_\beta \, f_{\delta - \alpha - \beta}, && f_{\delta - \beta} = f_{\delta - \alpha - \beta} \, f_\alpha - q \, f_\alpha \, f_{\delta - \alpha - \beta}. 
\end{align*}
The root vectors corresponding to the roots $\pm \delta$ are indexed by the elements of $\triangle_+$ and defined by the relations\footnote{The prime stands to distinguish two types of the root vectors, see the definition below.}
\begin{equation*}
e'_{\delta, \, \gamma} = e_\gamma e_{\delta - \gamma} - q^{-2} e_{\delta - \gamma} e_\gamma, \qquad f'_{\delta, \, \gamma} = f_{\delta - \gamma} f_\gamma - q^{-2} f_\gamma f_{\delta - \gamma}.
\end{equation*}
Now we define the root vectors corresponding to the remaining roots $\pm (\gamma + k \delta)$ and $\pm ((\delta - \gamma) + k \delta)$ as follows
\begin{align}
&e_{\gamma + k \delta} = [2]_q^{-1} (e_{\gamma + (k - 1) \delta} \, e'_{\delta, \, \gamma} -  e'_{\delta, \, \gamma} \, e_{\gamma + (k - 1) \delta}), \label{ekd1} \\
&e_{(\delta - \gamma) + k \delta} = [2]_q^{-1} (e'_{\delta, \, \gamma} \, e_{(\delta - \gamma) + (k - 1) \delta} - e_{(\delta - \gamma) + (k - 1) \delta} \,  e'_{\delta, \, \gamma}), \label{ekd2} \\
&f_{\gamma + k \delta} = [2]_q^{-1} (f'_{\delta, \, \gamma} \, f_{\gamma + (k - 1) \delta} -  f_{\gamma + (k - 1) \delta} \, f'_{\delta, \, \gamma}), \\
&f_{(\delta - \gamma) + k \delta} = [2]_q^{-1} (f_{(\delta - \gamma) + (k - 1) \delta} \, f'_{\delta, \, \gamma} - f'_{\delta, \, \gamma} \, f_{(\delta - \gamma) + (k - 1) \delta}).
\end{align}
The last step is to define the root vectors corresponding to the roots $k \delta$. For $k > 0$ they are defined as
\begin{equation}
e'_{k \delta, \, \gamma} = e_{\gamma + (k - 1) \delta} \, e_{\delta - \gamma} - q^{-2} e_{\gamma - \delta} \, e_{\gamma + (k - 1) \delta}, \label{epkd}
\end{equation}
and for $k < 0$ as
\begin{equation*}
f'_{k \delta, \, \gamma} = f_{\delta - \gamma} \, f_{\gamma + (k - 1) \delta} - q^2 f_{\gamma + (k - 1) \delta} \, f_{\gamma - \delta},
\end{equation*}
where $\gamma \in \triangle_+$. The second type of vectors corresponding to the roots $k \delta$ is defined by the equation
\begin{equation}
e_{\delta, \, \gamma}(\zeta) = \kappa_q^{-1} \log \bigl(1 + \kappa_q e'_{\delta, \, \gamma}(\zeta)\bigr), \qquad f_{\delta, \, \gamma}(\zeta) = - \kappa_q^{-1} \log \bigl(1 - \kappa_q f'_{\delta, \, \gamma}(\zeta)\bigr), \label{pwp}
\end{equation}
where
\begin{gather*}
e'_{\delta, \, \gamma}(\zeta) = \sum_{k = 1}^\infty e'_{k \delta, \, \gamma} \, \zeta^k, \qquad e_{\delta, \, \gamma}(\zeta) = \sum_{k = 1}^\infty e_{k \delta, \, \gamma} \, \zeta^k, \\
f'_{\delta, \, \gamma}(\zeta) = \sum_{k = 1}^\infty f'_{k \delta, \, \gamma} \, \zeta^k, \qquad f_{\delta, \, \gamma}(\zeta) = \sum_{k = 1}^\infty f_{k \delta, \, \gamma} \, \zeta^k.
\end{gather*}
It is useful to have in mind that all $e_{k \delta, \, \gamma}$ commute and all $f_{k \delta, \, \gamma}$ commute as well. 

The next ingredient of the Khoroshkin--Tolstoy construction is a normal order of $\triangle_+$. We use the following one \cite{BraGouZha95}
\begin{multline*}
\alpha, \ \alpha + \beta, \ \alpha + \delta, \alpha + \beta + \delta, \ \alpha + 2 \delta, \ \ \alpha + \beta + 2 \delta, \ \ldots \ , \ \beta, \ \beta + \delta, \ \beta + 2 \delta, \ \ldots \ , \\ 
\delta, \ 2 \delta, \ \ldots \ , \ \ldots \ , \ (\delta - \beta) + 2 \delta, \ (\delta - \beta) + \delta, \ \delta - \beta, \ \ldots \ , \\ (\delta - \alpha) + 2 \delta, \ (\delta - \alpha - \beta) + 2 \delta, \ (\delta - \alpha) + \delta, \ (\delta - \alpha - \beta) + \delta, \ \delta - \alpha, \ \delta - \alpha - \beta.
\end{multline*}

After all, $\calB$ is constructed as the product of three factors
\begin{equation}
\calB = \calR_{\prec \delta} \, \calR_{\sim \delta} \, \calR_{\succ \delta}. \label{brrrr}
\end{equation}
The factor $\calR_{\prec \delta}$ is the product over $\gamma \in \triangle_+$ and $k \in \bbZ_+$ of the $q$-exponentials
\begin{equation*}
\calR_{\gamma + k \delta} = \exp_{q^{-2}} (\kappa_q \, e_{\gamma + k \delta} \otimes f_{\gamma + k \delta}).
\end{equation*}
Here and below we assume that the $q$-exponential is defined as
\begin{equation*}
\exp_q(x) = \sum_{n = 0}^\infty q^{- (n - 1)n / 4} \frac{x^n}{[n]_q!} 
\end{equation*}
with
\begin{equation*}
[n]_q! = [n]_q [n - 1]_q \ldots [1]_q.
\end{equation*}
The order of the factors in $\calR_{\prec \delta}$ coincides with the chosen normal order of the roots $\gamma + k \delta$. For the second factor in (\ref{brrrr}) we have
\begin{equation}
\calR_{\sim \delta} = \exp \biggl(\kappa_q \sum_{k = 1}^\infty \sum_{i, j = 1}^2 u_{k, ij} \, e_{k \delta, \, \alpha_i} \otimes f_{k \delta, \alpha_j} \biggr), \label{rsd}
\end{equation}
where the quantities $u_{k, ij}$ are the entries of the matrix
\begin{equation}
u_k = \frac{k}{[k]_q} \frac{1}{q^{2k} + 1 + q^{-2k}} \left(\begin{array}{cc}
q^k + q^{-k} & (-1)^k \\
(-1)^k & q^k + q^{-k}
\end{array} \right). \label{uk}
\end{equation}
The last factor in (\ref{brrrr}) is the product over $\gamma \in \triangle_+$ and $k \in \bbZ_+$ of the $q$-exponentials
\begin{equation*}
\calR_{(\delta - \gamma) + k \delta} = \exp_{q^{-2}} \bigl(\kappa_q \, e_{(\delta - \gamma) + k \delta} \otimes f_{(\gamma - \delta) + k \delta}\bigr).
\end{equation*}
The order of the factors in $\calR_{\succ \delta}$ coincides with the chosen normal order of the roots $(\delta - \gamma) + k \delta$.

\section{\texorpdfstring{$R$-matrix and monodromy operators}{R-matrix and monodromy operators}}

We construct objects related to integrability by choosing representations for the factors of the tensor product $\uqlsliii \otimes \uqlsliii$ and then applying them to the universal $R$-matrix \cite{BooGoeKluNirRaz12a}. A spectral parameter is introduced by endowing $\uqlsliii$ with a $\bbZ$-gradation. Here we use the following procedure \cite{BooGoeKluNirRaz12a}. Given $\zeta \in \bbC^\times$, we define an automorphism $\Gamma_\zeta$ of $\uqlsliii$ by its action on the generators of $\uqlsliii$ as
\begin{equation}
\Gamma_\zeta(q^x) = q^x, \qquad \Gamma_\zeta(e_i) = \zeta^{s_i} e_i, \qquad \Gamma_\zeta(f_i) = \zeta^{-s_i} f_i, \label{phizeta}
\end{equation}
where $s_i$ are arbitrary integers. The automorphisms $\Gamma_\zeta$ corresponds to the $\bbZ$-gradation with the grading subspaces
\begin{equation*}
\uqlsliii_m = \{ a \in \uqlsliii \mid \Gamma_\zeta(a) = \zeta^m a \}.
\end{equation*}
Note that for any $\zeta \in \bbC^\times$ the universal $R$-matrix of $\uqlsliii$ satisfies the condition
\begin{equation}
(\Gamma_\zeta \otimes \Gamma_\zeta)(\calR) = \calR. \label{phiphir}
\end{equation}
Below we use the notation
\begin{equation*}
s_\delta = s_0 + s_1 + s_2, \qquad s_\alpha = s_1, \qquad s_\beta = s_2.
\end{equation*}

\subsection{\texorpdfstring{$R$-matrix}{R-matrix}} \label{s:rm}

The first useful object is an $R$-operator, or the $R$-matrix associated with it. To define the corresponding representation we use the Jimbo's homomorphism
\begin{equation*}
\varphi: \uqlsliii \to \uqgliii
\end{equation*}
defined by the relations
\begin{align}
& \varphi(q^{\nu h_{\delta - \alpha - \beta}}) = q^{\nu(G_3 - G_1)}, && \varphi(q^{\nu h_\alpha}) = q^{\nu (G_1 - G_2)}, && \varphi(q^{\nu h_\beta}) = q^{\nu (G_2 - G_3)}, \label{j2h} \\
& \varphi(e_{\delta - \alpha - \beta}) = F_3 \, q^{- G_1 - G_3}, && \varphi(e_\alpha) = E_1, && \varphi(e_\beta) = E_2, \label{j2e} \\
& \varphi(f_{\delta - \alpha - \beta}) = E_3 \, q^{G_1 + G_3} , && \varphi(f_\alpha) = F_1, && \varphi(f_\beta) = F_2, \label{j2f}
\end{align}
see the paper \cite{Jim86a}. Note that this is not a homomorphism of Hopf algebras.

We denote by $\widetilde \pi^\lambda$ the infinite-dimensional representation of $\uqgliii$ with the highest weight $\lambda = (\lambda_1, \lambda_2, \lambda_3)$, $\lambda_1, \lambda_2, \lambda_3 \in \bbC$. It is a highest weight representation with the highest weight vector $v_0$ characterised by the equalities
\begin{equation*}
q^{\nu G_1} v_0 = q^{\nu \lambda_1} v_0, \quad q^{\nu G_2} v_0 = q^{\nu \lambda_2} v_0, \quad q^{\nu G_3} v_0 = q^{\nu \lambda_3} v_0, \quad E_1 v_0 = 0, \quad E_2 v_0 = 0.
\end{equation*}
Using the Jimbo's homomorphism we define the infinite-dimensional representation
\begin{equation*}
\widetilde \varphi^\lambda = \widetilde \pi^\lambda \circ \varphi
\end{equation*}
of $\uqlsliii$. When $\lambda_1 - \lambda_2$ and $\lambda_2 - \lambda_3$ are non-negative integers the representation $\widetilde \pi^\lambda$ has the maximal subrepresentation such that the corresponding quotient representation is finite-dimensional. We denote this representation of $\uqgliii$ by $\pi^\lambda$ and define the corresponding representation of $\uqlsliii$ as
\begin{equation*}
\varphi^\lambda = \pi^\lambda \circ \varphi.
\end{equation*}

We consider the $R$-matrix defined by the representation $\pi^{(1,0,0)}$ of $\uqgliii$ that can be realised on the space $\bbC^3$ as
\begin{gather}
\pi^{(1,0,0)} (q^{\nu G_1}) = q^\nu E_{11} + E_{22} + E_{33}, \qquad \pi^{(1,0,0)} (q^{\nu G_2}) = E_{11} + q^\nu E_{22} + E_{33}, \label{pib} \\
\pi^{(1,0,0)} (q^{\nu G_3}) = E_{11} + E_{22} + q^\nu E_{33}, \\
\pi^{(1,0,0)}(E_1) = E_{12}, \quad \pi^{(1,0,0)}(E_2) = E_{23}, \quad
\pi^{(1,0,0)}(F_1) = E_{21}, \quad \pi^{(1,0,0)}(F_2) = E_{32}, \label{pie}
\end{gather}
where $E_{ab} \in \End(\bbC^3)$, $a, b = 1, 2, 3$, are defined by their action on the vectors of the standard basis $\{e_a\}$ of $\bbC^3$:
\begin{equation*}
E_{ab} e_c = \delta_{bc} e_a.
\end{equation*}
Using the well known relation
\begin{equation*}
E_{ab} E_{cd} = \delta_{bc} E_{ad},
\end{equation*}
we see that
\begin{equation}
\pi^{(1,0,0)}(E_3) = E_{13}, \qquad \pi^{(1,0,0)}(F_3) = E_{31}. \label{pi3}
\end{equation}
The $R$-operator associated with the representation $\pi^{(1,0,0)}$ is defined as\footnote{Remind that $\varphi^\lambda = \pi^\lambda \circ \varphi$.}
\begin{equation*}
R(\zeta_1 | \zeta_2) = \bigl((\varphi^{(1,0,0)} \circ \Gamma_{\zeta_1}) \otimes (\varphi^{(1,0,0)} \circ \Gamma_{\zeta_2}) \bigr)(\calR).
\end{equation*}
It follows from (\ref{phiphir}) that
\begin{equation*}
R(\zeta_1 \nu | \zeta_2 \nu) = R(\zeta_1 | \zeta_2),
\end{equation*}
therefore, one has
\begin{equation*}
R(\zeta_1 | \zeta_2) = R(\zeta_1^{\strut} \zeta_2^{-1}),
\end{equation*}
where $R(\zeta) = R(\zeta | 1)$. One can write
\begin{equation}
R(\zeta) = \sum_{a, b, c, d = 1}^3 E_{ac} \otimes E_{bd} \, \bm R_{ab | cd}(\zeta) \label{rme}
\end{equation}
and define the $9 \times 9$ matrix
\begin{equation*}
\bm R(\zeta) = (\bm R_{ab | cd}(\zeta)).
\end{equation*}
It can be shown \cite{BraGouZha95, BooGoeKluNirRaz10} that 
\begin{equation}
\bm R(\zeta) = \bm B(\zeta) \bm K. \label{rfrd}
\end{equation}
Here for the non-zero entries of the matrix $\bm B(\zeta) = (\bm B_{a b | c d}(\zeta))$  we have
\begin{align}
& \bm B_{11 | 11}(\zeta) = \bm B_{22 | 22}(\zeta) = \bm B_{33 | 33}(\zeta) = \rme^{f(\zeta^{s_\delta})} (1- q^{-2} \zeta^{s_\delta}), \label{brb} \\*
& \bm B_{12 | 12}(\zeta) = \bm B_{13 | 13}(\zeta) = \bm B_{21 | 21}(\zeta) \notag \\*
& \hspace{2.7em} {} = \bm B_{23 | 23}(\zeta) = \bm B_{31 | 31}(\zeta) = \bm B_{32 | 32}(\zeta) =  \rme^{f(\zeta^{s_\delta})} (1 - \zeta^{s_\delta}), \\*
& \bm B_{12 | 21}(\zeta) =  \rme^{f(\zeta^{s_\delta})} \zeta^{s_\alpha} \kappa_q, \qquad \hspace{2.5em} \bm B_{13 | 31}(\zeta) =  \rme^{f(\zeta^{s_\delta})} \zeta^{s_\alpha + s_\beta} \kappa_q, \\*
& \bm B_{23 | 32}(\zeta) =  \rme^{f(\zeta^{s_\delta})} \zeta^{s_\beta} \kappa_q, \qquad \hspace{2.5em} \bm B_{21 | 12}(\zeta) =  \rme^{f(\zeta^{s_\delta})} \zeta^{s_\delta - s_\alpha} \kappa_q, \\*
& \bm B_{31 | 13}(\zeta) =  \rme^{f(\zeta^{s_\delta})} \zeta^{s_\delta - s_\alpha - s_\beta} \kappa_q, \qquad \bm B_{32 | 23}(\zeta) =  \rme^{f(\zeta^{s_\delta})} \zeta^{s_\delta - s_\beta} \kappa_q, \label{bre}
\end{align}
where $f(\zeta)$ is a transcendental function having the form
\begin{equation*}
f(\zeta) = f_3(q^2 \zeta) + f_3(\zeta) + f_3(q^{-4} \zeta))
\end{equation*}
with
\begin{equation}
f_3(\zeta) = \sum_{k = 0}^\infty \frac{1}{q^{2k} + 1 + q^{-2k}} \frac{\zeta^k}{k} = \sum_{k = 0}^\infty \frac{1}{[3]_{q^k}} \frac{\zeta^k}{k}. \label{l3}
\end{equation}
The non-zero entries of the matrix $\bm K = (\bm K_{a b | c d})$ are
\begin{align}
&\bm K_{11 | 11} = \bm K_{22 | 22} = \bm K_{33 | 33} = q^{2/3}, \label{db} \\*
&\bm K_{12 | 12} = \bm K_{13 | 13} = \bm K_{21 | 21} = \bm K_{23 | 23} = \bm K_{31 | 31} = \bm K_{32 | 32} = q^{-1/3}. \label{de}
\end{align}

In fact, to define an $R$-operator one can use any finite-dimensional representation of $\uqgliii$. Different choices correspond to different quantum integrable models. In this work we restrict ourselves to the models defined by the representation $\pi^{(1,0,0)}$. 

\subsection{Ansatz for basic monodromy operator}

To construct monodromy operators one uses different representations for different factors of the tensor product $\uqlsliii \otimes \uqlsliii$. We use for the first factor an infinite-dimensional representation $\widetilde \varphi^\lambda$ or a finite-dimensional representation $\varphi^\lambda$, for the second factor the representation $\varphi^{(1,0,0)}$, and denote
\begin{align*}
\widetilde M^\lambda(\zeta | \eta) &= \bigl( (\widetilde \varphi^\lambda \circ \Gamma_\zeta) \otimes (\varphi^{(1,0,0)} \circ \Gamma_\eta) \bigr)(\calR), \\ M^\lambda(\zeta | \eta) &= \bigl( (\varphi^\lambda \circ \Gamma_\zeta)) \otimes (\varphi^{(1,0,0)} \circ \Gamma_\eta) \bigr)(\calR).
\end{align*}
From the point of view of spin chains such an operator corresponds to a one-site chain. In general, for a chain of length $n$ one takes instead of the representation $\varphi^{(1,0,0)} \circ \Gamma_\eta$ the tensor product of the representations $\varphi^{(1,0,0)} \circ \Gamma_{\eta_1}$, \ldots , $\varphi^{(1,0,0)} \circ \Gamma_{\eta_n}$, see, for example, the paper \cite{BooGoeKluNirRaz12a}. 

In fact, it is convenient for applications to define the monodromy operator
\begin{equation*}
M(\zeta | \eta) = ((\varphi \circ \Gamma_\zeta) \otimes (\varphi^{(1,0,0)} \circ \Gamma_\eta))(\calR),
\end{equation*}
and use the relations
\begin{equation*}
\widetilde M^\lambda(\zeta | \eta) = (\widetilde \pi^\lambda \otimes \id)(M(\zeta | \eta)), \qquad M^\lambda(\zeta | \eta) = (\pi^\lambda \otimes \id)(M(\zeta | \eta)).
\end{equation*}
Using equation (\ref{phiphir}), one can demonstrate that
\begin{equation*}
M(\zeta \nu | \eta \nu) = M(\zeta | \eta)
\end{equation*}
for any $\nu \in \bbC^\times$. Therefore, one has
\begin{equation*}
M(\zeta | \nu) = M(\zeta \eta^{-1}),
\end{equation*}
where $M(\zeta) = M(\zeta | 1)$. It is clear that
\begin{equation}
M(\zeta) = ((\varphi \circ \Gamma_\zeta) \otimes \varphi^{(1,0,0)})(\calR). \label{mzeta}
\end{equation}
We call the operator $M(\zeta)$ the basic monodromy operator corresponding to the automorphism $\varphi$. It can be represented as
\begin{equation}
M(\zeta) = \sum_{a, b = 1}^3 \bbM_{a b}(\zeta) \otimes E_{a b}, \label{mme}
\end{equation}
where $\bbM_{a b}(\zeta)$ are some unique elements of $\uqgliii$. Introduce the matrix
\begin{equation*}
\bbM(\zeta) = (\bbM_{ab}(\zeta)).
\end{equation*}
By definition, it is an element of $\Mat_3(\uqgliii)$. It is evident that $\bbM(\zeta)$ contains the same information as $M(\zeta)$. Slightly abusing terminology, we call $\bbM(\zeta)$ a basic monodromy operator. Let $\bm P = (\bm P_{ab | cd})$ be the matrix defined by the equality
\begin{equation*}
\bm P_{ab | cd} = \delta_{ad} \, \delta_{bc}.
\end{equation*}
The corresponding linear operator is the permutation of the factors of the tensor product $\bbC^3 \otimes \bbC^3$. It follows from the Yang--Baxter equation (\ref{rrr}) that
\begin{equation}
\hat {\bm R}(\zeta_1^{\mathstrut} \zeta_2^{-1}) \bigl( \bbM(\zeta_1) \boxtimes \bbM(\zeta_2) \bigr) = \bigl( \bbM(\zeta_2) \boxtimes \bbM(\zeta_1) \bigr) \hat {\bm R}(\zeta_1^{\mathstrut}  \zeta_2^{-1}), \label{rmm}
\end{equation}
where $\hat {\bm R}(\zeta) = \bm R(\zeta) \bm P$, and $\boxtimes$ is a natural generalisation of the Kronecker product to the case of matrices with noncommuting entries, see, for example, the paper \cite{BooGoeKluNirRaz12a}.

Using (\ref{mme}), we write
\begin{equation*}
M^\lambda(\zeta) = \sum_{a, b = 1}^3 \pi^\lambda(\bbM_{a b}(\zeta)) \otimes E_{a b} = \sum_{a, b = 1}^3 \bbM_{a b}^\lambda(\zeta) \otimes E_{a b},
\end{equation*}
and define the matrix
\begin{equation*}
\bbM^\lambda(\zeta) = (\bbM_{ab}^\lambda(\zeta)).
\end{equation*}
It is clear that
\begin{equation*}
R(\zeta) = M^{(1,0,0)}(\zeta) = (\pi^{(1,0,0)} \otimes \id)(M(\zeta)).
\end{equation*}
Therefore,
\begin{equation}
\bbM^{(1,0,0)}(\zeta) = \bbR(\zeta), \label{mer}
\end{equation}
where
\begin{equation*}
\bbR(\zeta) = (\bbR_{ab}(\zeta)),
\end{equation*}
and the quantities $\bbR_{ab}(\zeta) \in \Mat_3(\bbC)$ are defined by the equality
\begin{equation*}
R(\zeta) = \sum_{a,b = 1}^3 \bbR_{ab}(\zeta) \otimes E_{ab} = \sum_{a,b = 1}^3 \bigl( \sum_{c, d = 1}^3  \bm R_{ca | db}(\zeta) E_{cd} \bigr) \otimes E_{ab},
\end{equation*}
see (\ref{rme}) for the definition of $\bm R_{ab | cd}(\zeta)$. In accordance with (\ref{rfrd}) we have
\begin{equation*}
\bbR(\zeta) = \bbB(\zeta) \bbK.
\end{equation*}
Here $\bbB(\zeta) = (\bbB_{ab}(\zeta))$ and $\bbK = (\bbK_{ab})$, where, as follows from (\ref{brb})--(\ref{bre}), (\ref{pib})--(\ref{pie}) and (\ref{pi3}), the matrices $\bbB_{ab}(\zeta) \in \Mat_3(\bbC)$ can be represented as
\begin{align}
\bbB_{11}(\zeta) &= \rme^{f(\zeta^{s_\delta})} \pi^{(1,0,0)} \bigl(1 - \zeta^{s_\delta} q^{-2 G_1} \bigr), \label{obbrb} \\
\bbB_{22}(\zeta) &= \rme^{f(\zeta^{s_\delta})} \pi^{(1,0,0)} \bigl(1 - \zeta^{s_\delta} q^{- 2 G_2} \bigr), \\
\bbB_{33}(\zeta) &= \rme^{f(\zeta^{s_\delta})} \pi^{(1,0,0)} \bigl(1 - \zeta^{s_\delta} q^{- 2 G_3} \bigr), \\
\bbB_{12}(\zeta) &= \rme^{f(\zeta^{s_\delta})} \zeta^{s_\delta - s_\alpha} \kappa_q \, \pi^{(1,0,0)}(F_1), \quad \bbB_{13}(\zeta) = \rme^{f(\zeta^{s_\delta})} \zeta^{s_\delta - s_\alpha - s_\beta} \kappa_q \pi^{(1,0,0)}(F_3), \\
\bbB_{23}(\zeta) &= \rme^{f(\zeta^{s_\delta})} \zeta^{s_\delta - s_\beta} \kappa_q \, \pi^{(1,0,0)}(F_2), \quad \bbB_{21}(\zeta) = \rme^{f(\zeta^{s_\delta})} \zeta^{s_\alpha} \kappa_q \, \pi^{(1,0,0)}(E_1), \\
\bbB_{31}(\zeta) &= \rme^{f(\zeta^{s_\delta})} \zeta^{s_\alpha + s_\beta} \kappa_q \, \pi^{(1,0,0)}(E_3), \hspace{.8em} \bbB_{32}(\zeta) = \rme^{f(\zeta^{s_\delta})} \zeta^{s_\beta} \kappa_q \, \pi^{(1,0,0)}(E_2), \label{obbre}
\end{align}
while, as follows from (\ref{db}), (\ref{de}) and (\ref{pib}), for the matrices $\bbK_{ab} \in \Mat_3(\bbC)$ we have
\begin{gather}
\bbK_{11} = \pi^{(1,0,0)}(q^{-G/3 + G_1}), \qquad \bbK_{22} = \pi^{(1,0,0)}(q^{-G/3 + G_2}), \\
\bbK_{33} = \pi^{(1,0,0)}(q^{-G/3 + G_3}). \label{obbde}
\end{gather}
Equations (\ref{obbrb})--(\ref{obbde}) suggest to assume that
\begin{equation}
\bbM(\zeta) = \bbN(\zeta) \bbD. \label{afm}
\end{equation}
Here $\bbN(\zeta)$ is the matrix of the form
\begin{equation}
\bbN(\zeta) = \rme^{F(\zeta^{s_\delta})} \left( \begin{array}{rrr}
\bbN'_{11}(\zeta^{s_\delta}) & \zeta^{s_\delta - s_\alpha} \bbN'_{12} \phantom{(\zeta^{s_\delta})} & \zeta^{s_\delta - s_\alpha - s_\beta} \bbN'_{13} \phantom{(\zeta^{s_\delta})} \\[.5em]
\zeta^{s_\alpha} \bbN'_{21} \phantom{(\zeta^{s_\delta})} & \bbN'_{22}(\zeta^{s_\delta}) & \zeta^{s_\delta - s_\beta} \bbN'_{23} \phantom{(\zeta^{s_\delta})} \\[.5em]
\zeta^{s_\alpha + s_\beta} \bbN'_{31} \phantom{(\zeta^{s_\delta})} & \zeta^{s_\beta} \bbN'_{32} \phantom{(\zeta^{s_\delta})} & \bbN'_{33}(\zeta^{s_\delta})
\end{array} \right), \label{ma}
\end{equation}
where
\begin{gather}
\bbN'_{11}(\zeta) = 1 - \zeta q^{- 2 G_1}, \qquad \bbN'_{22}(\zeta) = 1 - \zeta q^{- 2 G_2}, \label{md1} \\
\bbN'_{33}(\zeta) = 1 - \zeta q^{- 2 G_3}, \label{md2} \\
\bbN'_{12} = c_1 \kappa_q \, F_1 \, q^{c_{11} G_1 + c_{12} G_2 + c_{13} G_3}, \qquad \bbN'_{21} = d_1 \kappa_q \, E_1 \, q^{d_{11} G_1 + d_{12} G_2 + d_{13} G_3}, \\
\bbN'_{23} = c_2 \kappa_q \, F_2 \, q^{c_{21} G_1 + c_{22} G_2 + c_{23} G_3}, \qquad \bbN'_{32} = d_2 \kappa_q \, E_2 \, q^{d_{21} G_1 + d_{22} G_2 + d_{23} G_3}, \\
\hspace{.1em} \bbN'_{13} = c_3 \kappa_q \, F_3 \, q^{c_{31} G_1 + c_{32} G_2 + c_{33} G_3}, \qquad \bbN'_{31} = d_3 \kappa_q \, \, E_3 \, q^{d_{31} G_1 + d_{32} G_2 + d_{33} G_3},
\end{gather}
$F(\zeta)$ belongs to the centre of $\uqgliii$, and $\bbD$ is a constant diagonal matrix with the diagonal entries
\begin{gather}
\bbD_{11} = q^{-G/3 + G_1}, \qquad \bbD_{22} = q^{-G/3 + G_2}, \qquad \bbD_{33} = q^{-G/3 + G_3}. \label{rbbkb}
\end{gather}
To find the numbers $c_i$, $d_j$, $c_{ij}$ and $d_{ij}$, we substitute the ansatz (\ref{afm}) into equation (\ref{rmm}). One can determine that, up to an automorphism of the form (\ref{aut}), this equality is satisfied if we put
\begin{align}
\bbN'_{12} &= \kappa_q \, q \, F_1 \, q^{- G_1 - G_2}, & \bbN'_{23} &= \kappa_q \, q \, F_2 \, q^{- G_2 - G_3}, & \bbN'_{13} &= \kappa_q \, q \, F_3 \, q^{- G_1 - G_3}, \label{robbmb} \\*
\bbN'_{21} &= \kappa_q \, E_1, & \bbN'_{32} &= \kappa_q \, E_2,  & \bbN'_{31} &= \kappa_q \, E_3. \label{robbme}
\end{align}
The quantity $F(\zeta)$ remains arbitrary. In the next section we prove that equations (\ref{rbbkb})--(\ref{robbme}) really describe the monodromy operator obtained from the universal $R$-matrix using the mapping $(\varphi \circ \Gamma_\zeta) \otimes \varphi^{(1,0,0)}$ and find the expression for $F(\zeta)$.

\subsection{Sketch of the proof} \label{s:sp}

The proof is rather cumbersome and quite technical, therefore we only describe the main steps and leave the details to the reader.

Let us first discuss the general structure of the monodromy operator $M(\zeta)$. We have
\begin{equation*}
M(\zeta) = U(\zeta) V(\zeta) W(\zeta) D,
\end{equation*}
where
\begin{align*}
& U(\zeta) = \bigl((\varphi \circ \Gamma_\zeta) \otimes \varphi^{(1,0,0)}\bigr)(\calR_{\prec \delta}), && V(\zeta) = \bigl((\varphi \circ \Gamma_\zeta) \otimes \varphi^{(1,0,0)}\bigr)(\calR_{\sim \delta}), \\
& W(\zeta) = \bigl((\varphi \circ \Gamma_\zeta) \otimes \varphi^{(1,0,0)}\bigr)(\calR_{\succ \delta}), && D = \bigl((\varphi \circ \Gamma_\zeta) \otimes \varphi^{(1,0,0)}\bigr)(\calK).
\end{align*}
Using relations similar to (\ref{mme}), we define the matrices $\bbU(\zeta)$, $\bbV(\zeta)$, $\bbW(\zeta)$ and $\bbD$ with the entries in $\uqsliii$ corresponding to the operators $U(\zeta)$, $V(\zeta)$, $W(\zeta)$ and $D$ respectively.

Using (\ref{kphipi}), (\ref{j2h}), (\ref{pib}) and taking into account that the endomorphism $E_{aa}$ is the projection on the vector $e_a$ of the standard basis of $\bbC^3$, we see that (\ref{rbbkb}) gives the right expression for the nonzero entries of the matrix $\bbD$. Thus, it remains to demonstrate that under the appropriate choice of $F(\zeta)$ the equality
\begin{equation}
\bbU(\zeta) \bbV(\zeta) \bbW(\zeta) = \bbN(\zeta), \label{uvw}
\end{equation}
where $\bbN(\zeta)$ is determined by equations (\ref{ma}) and (\ref{robbmb})--(\ref{robbme}), is true.

Note that for any $\gamma \in \triangle_+$ one has
\begin{equation*}
(\varphi \circ \Gamma_\zeta)(e_{\gamma + k \delta}) = \zeta^{s_\gamma + k s_\delta} \varphi(e_{\gamma + k \delta}), \qquad (\varphi \circ \Gamma_\zeta)(e_{(\delta - \gamma) + k \delta}) = \zeta^{(s_\delta - s_\gamma) + k s_\delta} \varphi(e_{(\delta - \gamma) + k \delta}).
\end{equation*}
Further, in the same way as in the paper \cite{BooGoeKluNirRaz10}, we obtain
\begin{align*}
& \varphi^{(1,0,0)}(f_{\alpha + k \delta}) = (-1)^k q^{2 k} E_{21}, && \varphi^{(1,0,0)}(f_{(\delta - \alpha) + k \delta}) = (-1)^k q^{2 k + 1} E_{12}, \\
& \varphi^{(1,0,0)}(f_{\beta + k \delta}) = q^{3 k} E_{32}, && \varphi^{(1,0,0)}(f_{(\delta - \beta) + k \delta}) = - q^{3 k + 2} E_{23}, \\
& \varphi^{(1,0,0)}(f_{\alpha + \beta + k \delta}) = (-1)^k q^{2 k} E_{31}, && \varphi^{(1,0,0)}(f_{(\delta - \alpha - \beta) + k \delta}) = (-1)^k q^{2 k + 1} E_{13}.
\end{align*}
For $\gamma \in \triangle_+$ denote
\begin{equation*}
\bbE_\gamma(\zeta) = \sum_{k = 0}^\infty \varphi(e_{\gamma + k \delta}) \, \zeta^k, \qquad \bbE_{\delta - \gamma}(\zeta) = \sum_{k = 0}^\infty \varphi(e_{(\delta - \gamma) + k \delta}) \, \zeta^k.
\end{equation*}
Using the properties of the endomorphisms $E_{ab}$, we see that the matrices $\bbU(\zeta)$ and $\bbW(\zeta)$ has the form
\begin{gather}
\bbU(\zeta) = \left( \begin{array}{rrl}
1 \hspace{3em} & 0 \hspace{3em} & \hspace{1em} 0 \\[.3em]
\zeta^{s_\alpha} \bbU'_{21}(\zeta^{s_\delta}) & 1 \hspace{3em} & \hspace{1em} 0 \\[.3em]
\zeta^{s_\alpha + s_\beta} \bbU'_{31}(\zeta^{s_\delta}) & \zeta^{s_\beta} \bbU'_{32}(\zeta^{s_\delta}) & \hspace{1em} 1 \hspace{1em}
\end{array} \right), \label{bbu} \\
\bbW(\zeta) = \left( \begin{array}{lrr}
\hspace{1em} 1 \hspace{3em} & \zeta^{s_\delta - s_\alpha} \bbW'_{12}(\zeta^{s_\delta}) & \zeta^{s_\delta - s_\alpha - s_\beta} \bbW'_{13}(\zeta^{s_\delta}) \\[.3em]
\hspace{1em} 0 & 1 \hspace{3em} & \zeta^{s_\delta - s_\beta} \bbW'_{23}(\zeta^{s_\delta}) \\[.3em]
\hspace{1em} 0 & 0 \hspace{3em} & 1 \hspace{3em}
\end{array} \right),  \label{bbw}
\end{gather}
where
\begin{align}
& \bbU'_{21}(\zeta) = \kappa_q \, \bbE_{\alpha}(- q^2 \zeta), && \bbW'_{12}(\zeta) = \kappa_q \, q \, \bbE_{\delta - \alpha}(- q^2 \zeta), \label{ueb} \\*
& \bbU'_{32}(\zeta) = \kappa_q \, \bbE_\beta(q^3 \zeta), && \bbW'_{23}(\zeta) = - \kappa_q \, q^2 \, \bbE_{\delta - \beta}(q^3 \zeta), \\*
& \bbU'_{31}(\zeta) = \kappa_q \, \bbE_{\alpha + \beta}(- q^2 \zeta), && \bbW'_{13}(\zeta) = \kappa_q \, q \, \bbE_{\delta - \alpha - \beta}(- q^2 \zeta). \label{uee}
\end{align}
Further, it is easy to get convinced that
\begin{equation*}
(\varphi \circ \Gamma_\zeta)(e_{k \delta, \, \gamma}) = \zeta^{k s_\delta} \varphi(e_{k \delta, \, \gamma})
\end{equation*}
for any $\gamma \in \triangle_+$. Similarly as in the paper \cite{BooGoeKluNirRaz10} we obtain
\begin{align}
& \varphi^{(1,0,0)}(f_{k \delta, \, \alpha}) = (-1)^{k - 1} \frac{[k]_q}{k} q^k (E_{11} - q^{2k} E_{22}), \label{fkda} \\
& \varphi^{(1,0,0)}(f_{k \delta, \, \beta}) = - \frac{[k]_q}{k} q^{2 k} (E_{22} - q^{2k} E_{33}). \label{fkdb}
\end{align}
Using these relations and the definition (\ref{rsd}), we see that the matrix $\bbV(\zeta)$ is of the diagonal form:
\begin{equation}
\bbV(\zeta) = \left( \begin{array}{ccc}
\bbV'_{11}(\zeta^{s_\delta}) & 0 & 0 \\[.3em]
0 & \bbV'_{22}(\zeta^{s_\delta}) & 0 \\[.3em]
0 & 0 & \bbV'_{33}(\zeta^{s_\delta}) \label{bbv}
\end{array} \right).
\end{equation}
while equation (\ref{uk}) gives
\begin{align}
\log \bbV'_{11}(\zeta) &= - \kappa_q \sum_{k = 1}^\infty \frac{\bigl((-1)^k (q^{4k} + q^{2k}) \varphi(e_{k \delta, \, \alpha}) +q^{3k} \varphi(e_{k \delta, \, \beta})\bigr) \, \zeta^k}{q^{4k} + q^{2k} + 1}, \label{lv1} \\
\log \bbV'_{22}(\zeta) &= \kappa_q \sum_{k = 1}^\infty \frac{\bigl((-1)^k q^{6k} \varphi(e_{k \delta, \, \alpha}) - q^{3k} \varphi(e_{k \delta, \, \beta}) \bigr) \, \zeta^k}{q^{4k} + q^{2k} + 1}, \label{lv2} \\
\log \bbV'_{33}(\zeta) &= \kappa_q \sum_{k = 1}^\infty \frac{\bigl((-1)^k q^{6k} \varphi(e_{k \delta, \, \alpha}) + (q^{7k} + q^{5k}) \varphi(e_{k \delta, \, \beta}) \bigr) \, \zeta^k}{q^{4k} + q^{2k} + 1}. \label{lv3}
\end{align}
It follows from (\ref{lv1}) and (\ref{lv2}) that
\begin{equation}
\log \bbV'_{22}(\zeta) - \log \bbV'_{11}(\zeta) = \kappa_q \bbE_{\delta, \, \alpha}(- q^2 \zeta). \label{lvlv}
\end{equation}
Here and below for $\gamma \in \triangle_+$ we use the notation
\begin{equation*}
\bbE_{\delta, \, \gamma}(\zeta) = \sum_{k = 1}^\infty \varphi(e_{k \delta, \, \gamma}) \, \zeta^k.
\end{equation*}
The definition (\ref{pwp}) gives
\begin{equation*}
1 + \kappa_q \bbE'_{\delta, \, \gamma}(\zeta) = \exp \bigl( \kappa_q \bbE_{\delta, \, \gamma}(\zeta) \bigr),
\end{equation*}
where
\begin{equation*}
\bbE'_{\delta, \, \gamma}(\zeta) = \sum_{k = 1}^\infty \varphi(e'_{k \delta, \, \gamma}) \, \zeta^k.
\end{equation*}
Hence, it follows from (\ref{lvlv}) that
\begin{equation}
1 + \kappa_q \bbE'_{\delta, \, \alpha}(\zeta) = \bbV'^{-1}_{11}(- q^{-2} \zeta) \bbV'^{\mathstrut}_{22}(- q^{-2} \zeta). \label{bea}
\end{equation}
In the same way, using (\ref{lv2}) and (\ref{lv3}), we obtain
\begin{equation}
1 + \kappa_q \bbE'_{\delta, \, \beta}(\zeta) = \bbV'^{-1}_{22}(q^{-3} \zeta) \bbV'^{\mathstrut}_{33}(q^{-3} \zeta). \label{beb}
\end{equation}
In fact, using different pairs of relations from (\ref{lv1})--(\ref{lv3}) we obtain different expressions for $\bbE'_{\delta, \, \gamma}(\zeta)$, $\gamma = \alpha, \beta$. However, they are related one to another by the identity
\begin{equation}
\bbV'_{11}(q^2 \zeta) \bbV'_{22}(\zeta) \bbV'_{33}(q^{-2} \zeta) = 1, \label{vvv}
\end{equation}
which can be obtained from (\ref{lv1})--(\ref{lv3}).

Rewrite equation (\ref{uvw}) in the component form and resolve the obtained equalities with respect to the entries of the matrices $\bbU$, $\bbV$ and $\bbW$. We come to the system
\begin{align}
& \bbU'_{21}(\zeta) = \bbN'_{21} \bbN'^{-1}_{11}(\zeta), && \bbU'_{31}(\zeta) = \bbN'_{31} \bbN'^{-1}_{11}(\zeta), \label{m1} \\
& \bbU'_{32}(\zeta) = \bbN''_{32}(\zeta) \bbN''^{-1}_{22}(\zeta), && \bbW'_{12}(\zeta) = \bbN'^{-1}_{11}(\zeta) \bbN'_{12}, \label{m2} \\
& \bbW'_{13}(\zeta) = \bbN'^{-1}_{11}(\zeta) \bbN'_{13}, && \bbW'_{23}(\zeta) = \bbN''^{-1}_{22}(\zeta) \bbN''_{23}(\zeta), \label{m3} \\
& \bbV'_{11}(\zeta) = \rme^{F(\zeta)} \bbN'_{11}(\zeta), && \bbV'_{22}(\zeta) = \rme^{F(\zeta)} \bbN''_{22}(\zeta), \label{m4} \\
& \bbV'_{33}(\zeta) = \rme^{F(\zeta)} \bbN'''_{33}(\zeta), \label{m5}
\end{align}
where
\begin{align*}
\bbN''_{23}(\zeta) &= \bbN'_{23} - \bbN'_{21} \bbN'^{-1}_{11}(\zeta) \bbN'_{13}, & \bbN''_{32}(\zeta) &= \bbN'_{32} - \zeta \, \bbN'_{31} \bbN_{11}^{\prime -1}(\zeta) \bbN'_{12}, \\
\bbN''_{22}(\zeta) &= \bbN'_{22}(\zeta) - \zeta \, \bbN'_{21} \bbN'^{-1}_{11}(\zeta) \bbN'_{12}, & \bbN''_{33}(\zeta) &= \bbN'_{33}(\zeta) - \zeta \, \bbN'_{31} \bbN'^{-1}_{11}(\zeta) \bbN'_{13},
\end{align*}
and
\begin{equation*}
\bbN'''_{33}(\zeta) = \bbN''_{33}(\zeta) - \zeta \, \bbN''_{32}(\zeta) \bbN''^{-1}_{22}(\zeta) \bbN''_{23}(\zeta).
\end{equation*}
It follows from (\ref{ueb})--(\ref{uee}) that equations (\ref{m1})--(\ref{m3}) are equivalent to the equalities
\begin{align}
& \bbE_\alpha(\zeta) = \kappa_q^{-1} \bbN'_{21} \bbN'^{-1}_{11}(- q^{-2} \zeta), \label{emb} \\*
& \bbE_\beta(\zeta) = \kappa_q^{-1} \bbN''_{32}(q^{-3} \zeta) \bbN''^{-1}_{22}(q^{-3} \zeta), \label{emm} \\*
& \bbE_{\alpha + \beta}(\zeta) = \kappa_q^{-1} \bbN'_{31} \bbN'^{-1}_{11}(- q^{-2} \zeta), \\
& \bbE_{\delta - \alpha}(\zeta) = \kappa_q^{-1} q^{-1} \bbN'^{-1}_{11}(- q^{-2} \zeta) \bbN'_{12}, \\*
& \bbE_{\delta - \beta}(\zeta) = - \kappa_q^{-1} q^{-2} \bbN''^{-1}_{22}(q^{-3} \zeta) \bbN''_{23}(q^{-3} \zeta), \\*
& \bbE_{\delta - \alpha - \beta}(\zeta) = \kappa_q^{-1} q^{-1} \bbN'^{-1}_{11}(- q^{-2} \zeta) \bbN'_{13}. \label{eme}
\end{align}
Equations (\ref{ekd1}) and (\ref{ekd2}) give
\begin{align*}
& \bbE_\gamma(\zeta) - \varphi(e_\gamma) = [2]_q^{-1} \zeta \, [\bbE_\gamma(\zeta), \, \varphi(e'_{\delta, \, \gamma})], \\
& \bbE_{\delta - \gamma}(\zeta) - \varphi(e_{\delta - \gamma}) = [2]_q^{-1} \zeta \, [\varphi(e'_{\delta, \, \gamma}), \, \bbE_{\delta - \gamma}(\zeta)].
\end{align*}
These relations determine $\bbE_\gamma(\zeta)$ and $\bbE_{\delta - \gamma}(\zeta)$ uniquely. One can verify that the right hand sides of (\ref{emb})--(\ref{eme}) satisfy them. Hence, the equalities (\ref{emb})--(\ref{eme}) are true and, therefore, equations (\ref{m1})--(\ref{m3}) are also true.

Consider now the first equality of (\ref{m4}). It is clear that $F(\zeta)$ can be represented as
\begin{equation}
F(\zeta) = \sum_{k = 1}^\infty \frac{F_k}{q^{2k} + 1 + q^{-2k}} \, \frac{\zeta^k}{k}, \label{lz}
\end{equation}
cf. the definition (\ref{l3}). Taking into account (\ref{md1}) and (\ref{lv1}), we see that the first equality of (\ref{m4}) is true if an only if
\begin{equation*}
F_k = (q^{2k} + 1 + q^{-2k}) q^{- 2 k G_1} - \kappa_q \, q^k k \bigl((-1)^k (q^k + q^{-k}) \varphi(e_{k \delta, \, \alpha}) + \varphi(e_{k \delta, \, \beta})\bigr).
\end{equation*}
It follows from the first relation of (\ref{pwp}) that
\begin{equation*}
e_{k \delta, \gamma} = \sum_{\ell_1 + 2 \ell_2 + \cdots + k \ell_k = k} \frac{(-\kappa_q)^{\ell_1 + \ell_2 + \cdots + \ell_k - 1}(\ell_1 + \ell_2 + \cdots + \ell_k - 1)!}{\ell_1! \ell_2! \ldots \ell_k!} \, e_{\delta,\gamma}^{\prime \ell_1} \, e_{2 \delta,\gamma}^{\prime \ell_2} \, \ldots \, e_{k \delta,\gamma}^{\prime \ell_k}.
\end{equation*}
In particular,
\begin{align*}
& e_{\delta, \gamma} = e'_{\delta, \gamma}, \\
& e_{2 \delta, \gamma} = e'_{2 \delta, \gamma} - \kappa_q (e'_{\delta, \gamma})^2 / 2, \\
& e_{3 \delta, \gamma} = e'_{3 \delta, \gamma} - \kappa_q e'_{\delta, \gamma} e'_{2 \delta, \gamma} + \kappa_q^2  (e'_{\delta, \gamma})^3 / 3, \\
& e_{4 \delta, \gamma} = e'_{4 \delta, \gamma} - \kappa_q e'_{\delta, \gamma} e'_{3 \delta, \gamma} - \kappa_q (e'_{2 \delta, \gamma})^2 + \kappa_q^2  (e'_{\delta, \gamma})^2 e'_{2 \delta, \gamma} - \kappa_q^3  (e'_{\delta, \gamma})^4 / 4.
\end{align*}
We use these equalities to calculate $F_k$ for small $k$. For $k = 1$ we obtain
\begin{equation*}
F_1 = C^{(1)},
\end{equation*}
where
\begin{multline}
C^{(1)} = q^{-2 G_1 - 2} + q^{-2 G_2} + q^{- 2 G_3 + 2} + \kappa_q^2 F_1 E_1 q^{- G_1 - G_2 - 1} \\
+ \kappa_q^2 F_2 E_2 q^{- G_2 - G_3 + 1} + \kappa_q^2 F_3 E_3 q^{- G_1 - G_3 + 1} - \kappa_q^3 F_3 E_1 E_2 q^{- G_1 - G_3}. \label{c1}
\end{multline}
Note that $C^{(1)}$ belongs to the centre of $\uqgliii$. For $k = 2$ we see that
\begin{equation*}
F_2 = 2 C^{(2)} + C^{(1)2},
\end{equation*}
where
\begin{multline}
C^{(2)} = - q^{- 2 G_1 - 2 G_2 - 2} - q^{- 2 G_1 - 2 G_3} - q^{- 2 G_2 - 2 G_3 + 2} - \kappa_q^2 F_1 E_1 q^{- G_1 - G_2 - 2 G_3 + 1} \\
- \kappa_q^2 F_2 E_2 q^{- 2 G_1 - G_2 - G_3 - 1} - \kappa_q^2 F_3 E_3 q^{- G_1 - 2 G_2 - G_3 + 1} - \kappa_q^3 F_1 F_2 E_3 q^{- G_1 - 2 G_2 - G_3 + 1}. \label{c2}
\end{multline}
The element $C^{(2)}$ also belongs to the centre of $\uqgliii$. Further calculations give
\begin{equation*}
F_3 = 3 C^{(3)} + 3 C^{(2)} C^{(1)} + C^{(1)3}, \qquad F_4 = 4 C^{(3)} C^{(1)} + 2 C^{(2)2} + 4 C^{(2)} C^{(1)2} + C^{(1)4},
\end{equation*}
where
\begin{equation}
C^{(3)} = q^{-2(G_1 + G_2 + G_3)}. \label{c3}
\end{equation}
It is natural to assume now that all $F_k$ are determined by the equality
\begin{equation}
\sum_{k = 1}^\infty F_k \frac{\zeta^k}{k} = - \log (1 - C^{(1)} \zeta - C^{(2)} \zeta^2 - C^{(3)} \zeta^3). \label{ck}
\end{equation}
Note that in this case $F(\zeta)$ is uniquely determined by the expansion (\ref{lz}) and by the relation
\begin{equation}
F(q^2 \zeta) + F(\zeta) + F(q^{-2} \zeta) = - \log (1 - C^{(1)} \zeta - C^{(2)} \zeta^2 - C^{(3)} \zeta^3). \label{fffc}
\end{equation}
Let us show that the above assumption allows to prove the validity of equations (\ref{m4}) and (\ref{m5}).

Note that the definition (\ref{epkd}) is equivalent to the equality
\begin{equation*}
\bbE'_{\delta, \, \gamma}(\zeta) = \zeta \bigl( \bbE_\gamma(\zeta) \varphi(e_{\delta - \gamma}) - q^{-2} \varphi(e_{\delta - \gamma}) \bbE_\gamma(\zeta) \bigr).
\end{equation*}
Using this relation, (\ref{emb}) and the equality
\begin{equation*}
\varphi(e_{\delta - \alpha}) = F_1 q^{- G_1 - G_2},
\end{equation*}
we find that
\begin{equation}
1 + \kappa \bbE'_{\delta, \, \alpha}(\zeta) = \bbN'^{-1}_{11}(-q^{- 2} \zeta) \bbN''_{22}(-q^{- 2} \zeta). \label{opeka}
\end{equation}
Comparing with (\ref{bea}), we see that
\begin{equation}
\bbV'^{-1}_{11}(\zeta) \bbV'_{22}(\zeta) = \bbN'^{-1}_{11}(\zeta) \bbN''_{22}(\zeta). \label{v1v2}
\end{equation}
In a similar way, using (\ref{emm}) and the equality
\begin{equation*}
\varphi(e_{\delta - \beta}) = - q F_2 q^{- G_2 - G_3} + \kappa_q q^{-1} F_3 E_1 q^{- G_1 - G_3},
\end{equation*}
we conclude that
\begin{equation}
1 + \kappa \bbE'_{\delta, \, \beta}(\zeta) = \bbN''^{-1}_{22}(q^{- 3} \zeta) \bbN'''_{33}(q^{- 3} \zeta). \label{opekb}
\end{equation}
Comparing with (\ref{beb}), we obtain
\begin{equation}
\bbV'^{-1}_{22}(\zeta) \bbV'_{33}(\zeta) = \bbN''^{-1}_{22}(\zeta) \bbN'''_{33}(\zeta). \label{v2v3}
\end{equation}
Note also that it follows from (\ref{opeka}) and (\ref{opekb}) that
\begin{equation}
[\bbN'_{11}(\zeta), \bbE'_{\delta, \gamma}(\zeta')] = 0, \qquad [\bbN''_{22}(\zeta), \bbE'_{\delta, \gamma}(\zeta')] = 0 \label{NE}
\end{equation}
for $\gamma = \alpha, \beta$.

Introduce the quantity $F'(\zeta)$, such that
\begin{equation}
\rme^{F'(\zeta)} = \bbV'_{11}(\zeta) \bbN'^{-1}_{11}(\zeta). \label{psi1}
\end{equation}
Then, it follows from (\ref{v1v2}) and (\ref{v2v3}) that
\begin{equation}
\rme^{F'(\zeta)} = \bbV'_{22}(\zeta) \label{psi2} \bbN''^{-1}_{22}(\zeta)
\end{equation}
and
\begin{equation}
\rme^{F'(\zeta)} = \bbV'_{33}(\zeta) \bbN'''^{-1}_{33}(\zeta). \label{psi3}
\end{equation}
One can demonstrate that
\begin{equation}
(1 - C^{(1)} \zeta - C^{(2)} \zeta^2 - C^{(3)} \zeta^3) = \bbN'_{11}(q^2 \zeta) \, \bbN''_{22}(\zeta) \, \bbN'''_{33}(q^{-2} \zeta). \label{omc}
\end{equation}
This relation, together with (\ref{psi1})--(\ref{psi3}) and (\ref{vvv}), gives
\begin{equation*}
F'(q^2 \zeta) + F'(\zeta) + F'(q^{-2} \zeta) = - \log (1 - C^{(1)} \zeta - C^{(2)} \zeta^2 - C^{(3)} \zeta^3).
\end{equation*}
Here we used the equalities
\begin{gather*}
[\bbN'_{11}(\zeta) \bbV'_{22}(\zeta')] = 0, \qquad [\bbN'_{11}(\zeta) \bbV'_{33}(\zeta')] = 0, \\
[\bbN''_{22}(\zeta) \bbV'_{22}(\zeta')] = 0,
\end{gather*}
which follow from (\ref{NE}). Thus, $F'(\zeta) = F(\zeta)$ and the equalities (\ref{m4}) and (\ref{m5}) are true.

The final result of our consideration is
\begin{equation*}
\bbM(\zeta) = q^{-G/3} \rme^{F(\zeta^{s_\delta})} \left( \begin{array}{ccc}
q^{G_1} - \zeta^{s_\delta} q^{-G_1} & \zeta^{s_\delta - s_\alpha} \kappa_q q^{-G_1} F_1 & \zeta^{s_\delta - s_\alpha - s_\beta} \kappa_q q^{-G_1} F_3 \\[.5em]
\zeta^{s_\alpha} \kappa_q E_1 q^{G_1} & q^{G_2} - \zeta^{s_\delta} q^{-G_2} & \zeta^{s_\delta - s_\beta} \kappa_q q^{-G_2} F_2 \\[.5em]
\zeta^{s_\alpha + s_\beta} \kappa_q E_3 q^{G_1} & \zeta^{s_\beta} \kappa_q E_2 q^{G_2} & q^{G_3} - \zeta^{s_\delta} q^{-G_3}
\end{array} \right).
\end{equation*}
Here $F(\zeta)$ has the form (\ref{lz}) where $F_k$ are determined by equation (\ref{ck}) with $C^{(1)}$, $C^{(2)}$ and $C^{(3)}$ given by (\ref{c1}), (\ref{c2}) and (\ref{c3}).

The matrices $\bbM^\lambda(\zeta)$ are obtained from $\bbM(\zeta)$ by applying to its matrix elements the mapping $\pi^\lambda$. It follows from (\ref{c1}), (\ref{c2}) and (\ref{c3}) that
\begin{equation*}
\pi^\lambda(1 - C^{(1)} \zeta - C^{(2)} \zeta^2 - C^{(3)} \zeta^3) = (1 - q^{- 2(\lambda_1 + 1)} \zeta)(1 - q^{- 2 \lambda_2} \zeta)(1 - q^{- 2(\lambda_3 -1)} \zeta).
\end{equation*}
Then equation (\ref{ck}) gives
\begin{equation*}
\pi^\lambda(F_k) = q^{- 2(\lambda_1 + 1) k} + q^{- 2 \lambda_2 k} + q^{- 2(\lambda_3 - 1) k}
\end{equation*}
and we come to the relation
\begin{equation*}
\pi^\lambda(F(\zeta)) = f_3(q^{- 2 (\lambda_1 + 1)} \zeta) + f_3(q^{- 2 \lambda_2} \zeta) + f_3(q^{- 2 (\lambda_3 - 1)} \zeta).
\end{equation*}
In particular, we have
\begin{equation*}
\pi^{(1, 0, 0)}(F(\zeta)) \\= f_3(q^{- 4} \zeta) + f_3(\zeta) + f_3(q^2 \zeta).
\end{equation*}
Using this equality, we can check the validity of equation (\ref{mer}).

Applying to the matrix elements of $\bbM(\zeta)$ the automorphism (\ref{aut}) with
\begin{gather*}
\nu_1 = q^{-1/2}, \qquad \nu_2 = q^{-1/2}, \\
\nu_{11} = -1/2, \quad \nu_{12} = 1/2, \quad \nu_{13} = 0, \quad \nu_{21} = 0, \quad \nu_{22} = -1/2, \quad \nu_{23} = 1/2,
\end{gather*}
we obtain
\begin{multline*}
\bbM(\zeta) = q^{-G/3} \rme^{F(\zeta^{s_\delta})} \\*
\times \left( \begin{array}{ccc}
q^{G_1} - \zeta^{s_\delta} q^{-G_1} & \zeta^{s_\delta - s_\alpha} \kappa_q q^{-(G_1 + G_2 - 1)/2} F_1 & \zeta^{s_\delta - s_\alpha - s_\beta} \kappa_q q^{-(G_1 + G_3 - 1)/2} F_3 \\[.5em]
\zeta^{s_\alpha} \kappa_q E_1 q^{(G_2 + G_1 - 1)/2} & q^{G_2} - \zeta^{s_\delta} q^{-G_2} & \zeta^{s_\delta - s_\beta} \kappa_q q^{-(G_2 + G_3 - 1)/2} F_2 \\[.5em]
\zeta^{s_\alpha + s_\beta} \kappa_q E_3 q^{(G_3 + G_1 - 1)/2} & \zeta^{s_\beta} \kappa_q E_2 q^{(G_3 + G_2 - 1)/2} & q^{G_3} - \zeta^{s_\delta} q^{-G_3}
\end{array} \right).
\end{multline*}
This expression is fully consistent with the formula given by Jimbo \cite{Jim86a}.

\section{More monodromy operators}

There are two special automorphisms of $\uqlsliii$. The first one is defined by the relations
\begin{align}
& \tau(e_{\delta - \alpha - \beta}) = e_{\delta - \alpha - \beta}, && \tau(e_\alpha) = e_\beta, && \tau(e_\beta) = e_\alpha, \label{taue} \\
& \tau(f_{\delta - \alpha - \beta}) = f_{\delta - \alpha - \beta}, && \tau(f_\alpha) = f_\beta, && \tau(f_\beta) = f_\alpha, \\
& \tau(q^{\nu h_{\delta - \alpha - \beta}}) = q^{\nu h_{\delta - \alpha - \beta}}, && \tau(q^{\nu h_\alpha}) = q^{\nu h_\beta}, && \tau(q^{\nu h_\beta}) = q^{\nu h_\alpha}, \label{tauh}
\end{align}
and the second one is given by
\begin{align}
& \sigma(e_{\delta - \alpha - \beta}) = e_\alpha, && \sigma(e_\alpha) = e_\beta, && \sigma(e_\beta) = e_{\delta - \alpha - \beta}, \label{sigmae} \\*
& \sigma(f_{\delta - \alpha - \beta}) = f_\alpha, && \sigma(f_\alpha) = f_\beta, && \sigma(f_\beta) = f_{\delta - \alpha - \beta}, \\*
& \sigma(q^{\nu h_{\delta - \alpha - \beta}}) = q^{\nu h_\alpha}, && \sigma(q^{\nu h_\alpha}) = q^{\nu h_\beta}, && \sigma(q^{\nu h_\beta}) = q^{\nu h_{\delta - \alpha - \beta}}. \label{sigmah}
\end{align}
One can use these automorphism to define additional monodromy operators.

We start with the automorphism $\tau$. Denote
\begin{equation*}
\ovarphi = \varphi \circ \tau
\end{equation*}
and define the following representations of $\uqlsliii$
\begin{equation*}
\widetilde {\ovarphi}{}^\lambda = \widetilde \pi^\lambda \circ \ovarphi, \qquad \ovarphi{}^\lambda = \pi^\lambda \circ \ovarphi
\end{equation*}
and the corresponding monodromy operators
\begin{equation*}
\widetilde{\oM}{}^\lambda(\zeta) = ((\widetilde{\ovarphi}{}^\lambda \circ \Gamma_\zeta) \otimes \varphi^{(1,0,0)})(\calR), \qquad \oM{}^\lambda(\zeta) = ((\ovarphi{}^\lambda \circ \Gamma_\zeta) \otimes \varphi^{(1,0,0)})(\calR).
\end{equation*}
As above, it is convenient to introduce the basic monodromy operator corresponding to the automorphism $\ovarphi$:
\begin{equation*}
\oM(\zeta) = ((\ovarphi \circ \Gamma_\zeta) \otimes \varphi^{(1,0,0)})(\calR),
\end{equation*}
cf. (\ref{mzeta}), and use it to construct the monodromy operators $\widetilde{\oM}{}^\lambda(\zeta)$ and $\oM{}^\lambda(\zeta)$.

It follows from (\ref{dqx}), (\ref{defi}) that
\begin{equation*}
(\tau \otimes \tau) \circ \Delta = \Delta \circ \tau.
\end{equation*}
Similarly, (\ref{doqx}) and (\ref{doefi}) give
\begin{equation*}
(\tau \otimes \tau) \circ \Delta^{\mathrm{op}} = \Delta^{\mathrm{op}} \circ \tau.
\end{equation*}
Using the definition of the universal $R$-matrix (\ref{urm}), we obtain the equality
\begin{equation*}
((\tau \otimes \tau)(\calR)) \Delta(\tau(a)) ((\tau \otimes \tau)(\calR))^{-1} = \Delta^{\mathrm{op}}(\tau(a)).
\end{equation*}
Taking into account the uniqueness theorem for the universal $R$-matrix \cite{KhoTol92}, we conclude that
\begin{equation*}
(\tau \otimes \tau)(\calR) = \calR.
\end{equation*}
This identity allows to rewrite the definition of $\oM(\zeta)$ in the form
\begin{equation*}
\oM(\zeta) = \bigl( (\varphi \circ (\tau \circ \Gamma_\zeta \circ \tau^{-1})) \otimes \ovarphi^{(1,0,0)} \bigr)(\calR).
\end{equation*}
We see that to construct $\oM(\zeta)$ one can use for the first factor of the tensor product $\uqlsliii \otimes \uqlsliii$ the same mapping $\varphi$, which was used for the construction of $M(\zeta)$ and for the second factor the representation $\ovarphi^{(1,0,0)}$. Here the integers $s_\alpha$ and $s_\beta$ should be interchanged. In fact, the above equality gives us a possibility to use for the construction of $\oM(\zeta)$ formulas of section \ref{s:sp}.

The monodromy operator $\oM(\zeta)$ can be represented as
\begin{equation*}
\oM(\zeta) = \oU(\zeta) \oV(\zeta) \oW(\zeta) \oD,
\end{equation*}
where
\begin{align*}
& \oU(\zeta) = \bigl( (\varphi \circ (\tau \circ \Gamma_\zeta \circ \tau^{-1})) \otimes \ovarphi^{(1,0,0)} \bigr)(\calR_{\prec \delta}), \\
& \oV(\zeta) = \bigl( (\varphi \circ (\tau \circ \Gamma_\zeta \circ \tau^{-1})) \otimes \ovarphi^{(1,0,0)} \bigr)(\calR_{\sim \delta}), \\
& \oW(\zeta) = \bigl( (\varphi \circ (\tau \circ \Gamma_\zeta \circ \tau^{-1})) \otimes \ovarphi^{(1,0,0)} \bigr)(\calR_{\succ \delta}), \\
& \oD = \bigl( (\varphi \circ (\tau \circ \Gamma_\zeta \circ \tau^{-1})) \otimes \ovarphi^{(1,0,0)} \bigr)(\calK).
\end{align*}
Using relations similar to (\ref{mme}), we introduce the matrix $\obbM(\zeta)$ containing the same information as $\oM(\zeta)$ and write
\begin{equation*}
\obbM(\zeta) = \obbU(\zeta) \obbV(\zeta) \obbW(\zeta) \obbD,
\end{equation*}
where $\obbU(\zeta)$, $\obbV(\zeta)$, $\obbW(\zeta)$ and $\obbD$ are the analogues of $\bbU(\zeta)$, $\bbV(\zeta)$, $\bbW(\zeta)$ and $\bbD$ introduced in section \ref{s:sp}. It is clear that there are formulas similar to (\ref{bbu}), (\ref{bbw}) and (\ref{bbv}) with barred quantities.

One can demonstrate that
\begin{align*}
& \ovarphi^{(1,0,0)}(f_{k \delta, \, \alpha}) = - \frac{[k]_q}{k} q^{2 k} (E_{22} - q^{2k} E_{33}), \\
& \ovarphi^{(1,0,0)}(f_{k \delta, \, \beta}) = (-1)^{k - 1} \frac{[k]_q}{k} q^k (E_{11} - q^{2k} E_{22}).
\end{align*}
Comparing this equalities with (\ref{fkda}) and (\ref{fkdb}), we see that the expression for $\obbV(\zeta)$ can be obtained from the expression for $\bbV(\zeta)$ by interchanging $e_{k \delta, \alpha}$ and $e_{k \delta, \beta}$. In particular, we have
 \begin{equation*}
\obbV'_{11}(\zeta) = \bbV^{\prime -1}_{33}(- q^{-3} \zeta), \qquad \obbV'_{33}(\zeta) = \bbV^{\prime -1}_{11}(- q^3 \zeta).
\end{equation*}
Taking into account equations (\ref{m4}), (\ref{m5}), (\ref{omc}) and (\ref{fffc}), we obtain
\begin{align*}
& \obbV'_{11}(\zeta) = \rme^{F(- q^{-1} \zeta) + F(- q \zeta)} \, \bbN'_{11}(- q \zeta) \bbN''_{22}(- q^{-1} \zeta), \\
& \obbV'_{33}(\zeta) = \rme^{F(- q^{-1} \zeta) + F(- q \zeta)} \, \bbN''_{22}(- q \zeta) \bbN'''_{33}(- q^{-1} \zeta).
\end{align*}
It is not difficult to get convinced that
\begin{equation*}
\obbV'_{11}(q^2 \zeta) \obbV'_{22}(\zeta) \obbV'_{33}(q^{-2} \zeta) = 1.
\end{equation*}
Now equations (\ref{omc}) and (\ref{fffc}) give
\begin{equation*}
\obbV'_{22}(\zeta) = \rme^{F(- q^{-1} \zeta) + F(- q \zeta)} \, \bbN'_{11}(- q \zeta) \bbN'''_{33}(- q^{-1} \zeta).
\end{equation*}

Using the equalities
\begin{align*}
& \ovarphi^{(1,0,0)} (f_{\alpha + k \delta}) = q^{3k} E_{32}, && \ovarphi^{(1,0,0)} (f_{(\delta - \alpha) + k \delta}) = - q^{3k + 2} E_{23}, \\
& \ovarphi^{(1,0,0)} (f_{\beta + k \delta}) = (-1)^k q^{2k} E_{21}, && \ovarphi^{(1,0,0)} (f_{(\delta - \beta) + k \delta}) = (-1)^k q^{2k + 1} E_{12}, \\
& \ovarphi^{(1,0,0)} (f_{\alpha + \beta + k \delta}) = - q^{3k + 1} E_{31}, && \ovarphi^{(1,0,0)} (f_{(\delta - \alpha - \beta) + k \delta}) = q^{3k + 1} E_{13},
\end{align*}
and (\ref{emb})--(\ref{eme}), we obtain
\begin{align*}
& \obbU'_{21}(\zeta) = \bbN''_{32}(- q^{-1} \zeta) \bbN''^{-1}_{22}(- q^{-1} \zeta), \\
& \obbU'_{31}(\zeta) = - q \bbN'_{31} \bbN'^{-1}_{11}(- q \zeta) + \bbN'_{21} \bbN'^{-1}_{11}(-q \zeta) \bbN''_{32}(- q^{-1} \zeta) \bbN''^{-1}_{22}(- q^{-1} \zeta), \\
& \obbU'_{32}(\zeta) = \bbN'_{21} \bbN'^{-1}_{11}(-q \zeta), \\
& \obbW'_{12}(\zeta) = - q^{-1} \bbN''^{-1}_{22}(- q^{-1} \zeta) \bbN''_{23}(- q^{-1} \zeta), \\
& \obbW'_{13}(\zeta) = \bbN'^{-1}_{11}(-q \zeta) \bbN'_{13} + \zeta \bbN''^{-1}_{22}(- q^{-1} \zeta) \bbN''_{23}(- q^{-1} \zeta) \bbN'^{-1}_{11}(- q \zeta) \bbN'_{12}, \\
& \obbW'_{23}(\zeta) = - q \bbN'^{-1}_{11}(- q \zeta) \bbN'_{12}.
\end{align*}
Now we can determine that $\obbM(\zeta)$ has the form 
\begin{multline*}
\obbM(\zeta) = \rme^{\Psi(- q^{-1} \zeta^{s_\delta}) + \Psi(- q \zeta^{s_\delta})} q^{- 2 G / 3} \\* \times \left( \begin{array}{rrr}
\obbM'_{11}(\zeta^{s_\delta}) & \zeta^{s_\delta - s_\alpha} \obbM'_{12}(\zeta^{s_\delta}) & \zeta^{s_\delta - s_\alpha - s_\beta} \obbM'_{13}(\zeta^{s_\delta}) \\[.5em]
\zeta^{s_\alpha} \obbM'_{21}(\zeta^{s_\delta}) & \obbM'_{22}(\zeta^{s_\delta}) & \zeta^{s_\delta - s_\beta} \obbM'_{23} (\zeta^{s_\delta}) \\[.5em]
\zeta^{s_\alpha + s_\beta} \obbM'_{31}(\zeta^{s_\delta}) & \zeta^{s_\beta} \obbM'_{32}(\zeta^{s_\delta}) & \obbM'_{33}(\zeta^{s_\delta}),
\end{array} \right),
\end{multline*}
where
\begin{align*}
\obbM'_{11}(\zeta) &= (q^{G_1} + \zeta q^{- G_1 - 1}) (q^{G_2} + \zeta q^{- G_2 + 1}) + \zeta  \kappa_q^2  F_1 E_1, \\
\obbM'_{12}(\zeta) &= - \kappa_q q^2 (q^{G_1} + \zeta q^{- G_1 - 1}) F_2 q^{- G_2} + \kappa_q^2 F_3 E_1, \\
\obbM'_{13}(\zeta) &= \kappa_q q (q^{G_2} + \zeta q^{- G_2 + 1}) F_3 q^{- G_1} + \zeta \kappa_q^2 q^2 F_1 F_2 q^{- G_1 - G_2} , \\
\obbM'_{21}(\zeta) &= \kappa_q (q^{G_1} + \zeta q^{- G_1 - 1}) E_2 q^{G_2} + \zeta \kappa_q^2 F_1 E_3, \\
\obbM'_{22}(\zeta) &= (q^{G_1} + \zeta q^{- G_1 - 1}) (q^{G_3} + \zeta q^{- G_3 + 1}) + \zeta  \kappa_q^2  F_3 E_3, \\
\obbM'_{23}(\zeta) &= - \kappa_q (q^{G_3} + \zeta q^{- G_3 + 1}) F_1 q^{- G_1} + \kappa_q^2 q F_3 E_2 q^{- G_1 + G_2}, \\
\obbM'_{31}(\zeta) &= - \kappa_q q (q^{G_2} + \zeta q^{- G_2 - 1}) E_3 q^{G_1} + \kappa_q^2 E_1 E_2 q^{G_1 + G_2}, \\
\obbM'_{32}(\zeta) &= \kappa_q (q^{G_3} + \zeta q^{- G_3 + 1}) E_1 q^{G_1} + \zeta \kappa_q^2 q F_2 E_3 q^{G_1 - G_2}, \\
\obbM'_{33}(\zeta) &= (q^{G_2} + \zeta q^{- G_2 - 1}) (q^{G_3} + \zeta q^{- G_3 + 1}) + \zeta  \kappa_q^2  F_2 E_2.
\end{align*}

Returning to the automorphism $\sigma$, we define the homomorphisms
\begin{equation*}
\varphi_i = \varphi \circ \sigma^{- i + 1}, \qquad \ovarphi_i = \ovarphi \circ \sigma^{- i + 1},
\end{equation*}
where $i = 1, 2, 3$ and the corresponding basic monodromy operators
\begin{equation*}
M_i(\zeta) = \bigl( (\varphi_i \circ \Gamma_\zeta) \otimes \varphi^{(1,0,0)} \bigr)(\calR), \qquad \oM_i(\zeta) = \bigl( (\ovarphi_i \circ \Gamma_\zeta) \otimes \varphi^{(1,0,0)} \bigr)(\calR).
\end{equation*}
As in the case of the automorphism $\tau$, one can demonstrate that
\begin{equation*}
(\sigma \otimes \sigma)(\calR) = \calR.
\end{equation*}
Hence, one can write
\begin{equation*}
M_i(\zeta) = \bigl( (\varphi \circ (\sigma^{- i + 1} \circ \Gamma_\zeta \circ \sigma^{i - 1})) \otimes (\varphi^{(1,0,0)} \circ \sigma^{- i + 1}) \bigr)(\calR).
\end{equation*}
One has
\begin{equation*}
(\varphi^{(1,0,0)} \circ \sigma^{- i + 1})(a) = \Sigma^{-i + 1} \bigl( \varphi^{(1,0,0)}(a) \bigr) \Sigma^{i - 1},
\end{equation*}
where the matrix $\Sigma$ has the form
\begin{equation*}
\Sigma = \left( \begin{array}{ccc}
0 & 0 & 1 \\
1 & 0 & 0 \\
0 & 1 & 0                
\end{array} \right).
\end{equation*}
It is clear now that the monodromy operators $\bbM_i(\zeta)$ for different values of $i$ can be obtained by the similarity transformation generated by the corresponding power of $\Sigma$ and the corresponding power of the transformation
\begin{equation*}
s_\delta \to s_\delta, \qquad s_\alpha \to s_\beta, \qquad s_\beta \to s_\delta - s_\alpha - s_\beta.
\end{equation*}
The same is true for the monodromy operators $\obbM_i(\zeta)$.

\section{Conclusions}

Starting with the expression for the universal $R$-matrix given by Khoroshkin and Tolstoy \cite{TolKho92}, we constructed the basic monodromy operators for the case of the quantum group $\uqlsliii$. We see that despite of the fact that the formula given in \cite{TolKho92} is rather formal one can obtain explicit and sensible results. It is important that we have the exact result with the explicit form of the factors belonging to the centre of the quantum groups $\mathrm U_q(\mathfrak{sl}_3)$ and $\mathrm U_q(\mathfrak{gl}_3)$. An interesting by-product of our work is the expressions for quantum Casimir elements of $\mathrm U_q(\mathfrak{sl}_3)$ and $\mathrm U_q(\mathfrak{gl}_3)$.

{\em Acknowledgements.\/} This work was supported in part by the RFBR grants \#~10-01-00300 and \#~13-01-00217. The author would like to thank the Max Planck Institute for Mathematics in Bonn, where this work was started, for the hospitality extended to him during his stay there in February-May 2012.

\appendix

\stepcounter{section}

\section*{\texorpdfstring{Appendix. From $\uqgliii$ to $\uqsliii$}{Appendix. From Uq(gl3) to Uq(sl3)}}

To define basic monodromy operators, one can also use a homomorphism from the quantum group $\uqlsliii$ to the quantum group $\uqsliii$ defined by the relations
\begin{align}
& \varphi(q^{\nu h_0}) = q^{- \nu(H_1 + H_2)}, && \varphi(q^{\nu h_1}) = q^{\nu H_1}, && \varphi(q^{\nu h_2}) = q^{\nu H_2}, \label{slhb} \\
& \varphi(e_0) = F_3 \, q^{-(H_1 - H_2)/3}, && \varphi(e_1) = E_1, && \varphi(e_2) = E_2,  \\
& \varphi(f_0) = E_3 \, q^{(H_1 - H_2)/3} , && \varphi(f_1) = F_1, && \varphi(f_2) = F_2. \label{slhe}
\end{align}
Since we define the quantum group $\uqsliii$ as a subalgebra of $\uqgliii$, we can rewrite the expressions for $\varphi(e_0)$ and $\varphi(f_0)$ as
\begin{equation*}
\varphi(e_{\delta - \alpha - \beta}) = F_3 q^{- G_1 - G_3} q^{2 G / 3}, \qquad \varphi(f_{\delta - \alpha - \beta}) = E_3 q^{G_1 + G_3} q^{- 2 G / 3}.
\end{equation*}
Comparing these relations with the formulas describing the Jimbo's homomorphism, we conclude that the expressions for the basic monodromy operators based on the homomorphism defined by (\ref{slhb})--(\ref{slhe}) can be obtained from the expressions based on the Jimbo's homomorphism via the substitution
\begin{equation*}
\zeta^\delta \to \zeta^\delta q^{2 (G_1 + G_2 + G_3) / 3} q^{- 2 / 3}.
\end{equation*}
Now we see that the operator $\bbM(\zeta)$ has the form (\ref{afm}), where $\bbD$ is a diagonal matrix with the diagonal entries
\begin{equation*}
\bbD_{11} = q^{(2 H_1 + H_2)/3}, \qquad \bbD_{22} = q^{- (H_1 - H_2)/3}, \qquad \bbD_{33} = q^{-(H_1 + 2 H_2)/3}, \label{rbbk}
\end{equation*}
while the matrix $\bbN(\zeta)$ has the form (\ref{ma}), where
\begin{gather*}
\bbN'_{11}(\zeta) = 1 - \zeta q^{-(4 H_1 + 2 H_2 + 2) / 3}, \qquad \bbN'_{22}(\zeta) = 1 - \zeta q^{(2 H_1 - 2 H_2 - 2) / 3},\\
\bbN'_{33}(\zeta) = 1 - \zeta q^{(2 H_1 + 4 H_2 - 2) / 3},
\end{gather*}
and
\begin{align*}
\bbN'_{12} &= \kappa_q F_1 q^{- (H_1 + 2 H_2 - 1)/3}, && \bbN'_{21} = \kappa_q E_1, \\*
\bbN'_{23} &= \kappa_q F_2 q^{(2 H_1 + H_2 + 1)/3}, && \bbN'_{32} = \kappa_q E_2, \\*
\bbN'_{13} &= \kappa_q F_3 q^{- (H_1 - H_2 - 1)/3}, && \bbN'_{31} = \kappa_q E_3.
\end{align*}
For $F(\zeta)$ we again have the representation (\ref{lz}), where $F_k$ are determined by the equality (\ref{ck}) with
\begin{align*}
C^{(1)} &= q^{-(4H_1 + 2H_2 + 8)/3} + q^{(2H_1 - 2H_2 - 2)/3} + q^{(2H_1 + 4H_2 + 4)/3} \\*
& \hspace{5em} {} + \kappa_q^2 F_1 E_1 q^{-(H_1 + 2H_2 + 5)/3} + \kappa_q^2 F_2 E_2 q^{(2H_1 + H_2 + 1)/3} \\*
& \hspace{10em} {} + \kappa_q^2 F_3 E_3 q^{-(H_1 - H_2 - 1)/3} - \kappa_q^3 F_3 E_1 E_2 q^{-(H_1 - H_2 + 2)/3}, \\
C^{(2)} &= - q^{-(2H_1 + 4H_2 + 10)/3} - q^{-(2H_1 - 2H_2 + 4)/3} - q^{(4H_1 + 2H_2 + 2)/3} \\*
& \hspace{5em} {} - \kappa_q^2 F_1 E_1 q^{(H_1 + 2H_2 - 1)/3} - \kappa_q^2 F_2 E_2 q^{-(2H_1 + H_2 + 7)/3} \\*
& \hspace{10em} {} - \kappa_q^2 F_3 E_3 q^{(H_1 - H_2 - 1)/3} - \kappa_q^3 F_1 F_2 E_3 q^{(H_1 - H_2 - 1)/3}, \\
C^{(3)} &= q^{-2}.
\end{align*}
One can compare the obtained expressions for $C^{(1)}$ and $C^{(2)}$ with the expressions for the quantum Casimir operators from the paper \cite{Rod91}.

\providecommand{\href}[2]{#2}

\end{document}